\documentclass{article}

\usepackage{amscd,amsmath,amssymb,amsthm}
\usepackage{color}

\newcommand{\eps}{\varepsilon}

\newcommand{\R}{\mathbb R}
\newcommand{\N}{\mathbb N}

\newcommand{\then}{\Longrightarrow}
\newcommand{\J}{{\cal J}}

\DeclareMathOperator*{\esssup}{ess\; sup}

\DeclareMathOperator{\supp}{supp}
\newcommand\meas{{\rm meas}}
\newcommand\e{{\rm e}}

\newtheorem{corollary}{Corollary}[section]
\newtheorem{theorem}[corollary]{Theorem}
\newtheorem{lemma}[corollary]{Lemma}
\newtheorem{proposition}[corollary]{Proposition}

\theoremstyle{definition}
\newtheorem{definition}[corollary]{Definition}
\newtheorem{remark}[corollary]{Remark}
\newtheorem{example}[corollary]{Example}

\numberwithin{equation}{section}
\begin{document}

\title{{\bf A dichotomy result \\
for a modified Schr\"odinger equation\\
on unbounded domains} 
\footnote{The research that led to the present paper was partially supported 
by MUR--PRIN 2022 PNRR Project P2022YFAJH ``Linear and Nonlinear PDE's: 
New directions and Applications'' and
INdAM - GNAMPA Project 2024 ``Nonlinear problems in local and nonlocal settings with applications''.}}

\author{A.M. Candela, G. Palmieri and A. Salvatore\\
{\small Dipartimento di Matematica} \\
{\small Universit\`a degli Studi di Bari Aldo Moro} \\
{\small Via E. Orabona 4, 70125 Bari, Italy}\\
{\small \it annamaria.candela@uniba.it}\\
{\small \it giuliana.palmieri.uniba@gmail.com}\\ 
{\small \it addolorata.salvatore@uniba.it}}
\date{}

\maketitle

\begin{abstract}
This article aims to investigate the existence of  
bounded positive solutions of problem
\[ 
(P)\qquad 
\left\{
\begin{array}{ll}
- {\rm div} (a(x,u,\nabla u)) + A_t(x,u,\nabla u) = g(x,u)
  &\hbox{in $\Omega$,}\\
u\ = \ 0 & \hbox{on $\partial\Omega$,}
\end{array}
\right.
\]
with $A_t(x,t,\xi) = \frac{\partial A}{\partial t}(x,t,\xi)$, 
$a(x,t,\xi) = \nabla_\xi A(x,t,\xi)$
for a given $A(x,t,\xi)$ which grows as $|\xi|^p + |t|^p$ , $p > 1$,
where $\Omega \subseteq \R^N$, $N \ge 2$, is an open connected domain 
with Lipschitz boundary and infinite Lebesgue measure,
eventually $\Omega = \R^N$, which generalizes the modified Schr\"odinger equation
\[
- {\rm div} ((A^*_1(x) + A^*_2(x)|u|^{s}) \nabla u) + 
\frac{s}2 A^*_2(x)\ |u|^{s - 2} u\ |\nabla u|^2 + u\ 
=\ |u|^{\mu-2}u \quad\hbox{in $\R^3$.}
\]
Under suitable assumptions on $A(x,t,\xi)$ and $g(x,t)$, problem $(P)$  
has a variational structure. Then, even in lack of radial symmetry hypotheses,
one bounded positive solution of $(P)$ can be found by passing to 
the limit on a sequence $(u_k)_k$ of bounded solutions on bounded domains.
Furthermore, if stronger hypotheses are satisfied,
either such a solution is nontrivial or a constant $\bar{\lambda} > 0$ and a
sequence of points $(y_k)_k \subset \R^N$ exist such that 
\[
|y_k| \to +\infty\qquad \hbox{and}\qquad
\int_{B_1(y_k)} |u_k|^p dx \ge \bar{\lambda}\quad \hbox{for all $k \ge 1$.}
\]
\end{abstract}



\section{Introduction}

This article aims to investigate the existence of critical points of the nonlinear functional
\[
\J(u)\ =\ \int_\Omega A(x,u,\nabla u) dx - \int_\Omega G(x,u) dx,\qquad
u \in {\cal D}\subset W^{1,p}_0(\Omega),
\]
which generalizes the model problem 
\begin{equation}\label{modello}
\bar{\J}(u)\ =\ \frac1p\ \int_\Omega A_1(x,u) |\nabla u|^p dx
+ \frac1p\ \int_\Omega A_2(x,u) |u|^p dx - \int_\Omega G(x,u) dx,
\end{equation}
where $\Omega$ is an open connected domain in $\R^N$, $N \ge 2$, with Lipschitz boundary and infinite Lebesgue measure, 
eventually $\Omega = \R^N$, $p >1 $, and the functions
$A :\Omega \times \R \times \R^N \to \R$, respectively $A_j :\Omega \times \R \to \R$ with $j \in \{1,2\}$,
and $G : \Omega \times \R \to \R$ are given.

We note that, even in the simplest case when $\Omega$ is bounded
and $A_2(x,t)\equiv 0$, $G(x,t)\equiv 0$, with $A_1(x,t)$ smooth, the functional $\bar{\J}$
is not well defined in the whole space $W^{1,p}_0(\Omega)$
if $A_1(x,t)$ is unbounded with respect to $t$.
Moreover, even if $A_1(x,t)$ is bounded with respect to $t$
 but $\frac{\partial A_1}{\partial t}(x,t) \not\equiv 0$, 
then $\bar{\J}$ is defined in $W^{1,p}_0(\Omega)$ but  
is G\^ateaux differentiable only along directions of 
$W^{1,p}_0(\Omega) \cap L^\infty(\Omega)$. Thus, many authors have 
studied functional $\J$ by using nonsmooth techniques or introducing suitable changes of variables,
when possible, or giving a ``good'' definition of critical point either on bounded domains or in unbounded ones
(see, e.g., \cite{AB1, AG, BMP, BP, Ca, CD, CJ,CDM, LW, LWW, PSW, ShW, SC}).  

Here, as in \cite{CP2}, suitable assumptions assure that functional
$\J$ is $C^1$ in $X=W^{1,p}_0(\Omega) \cap L^\infty(\Omega)$ (see Proposition \ref{smooth1})
and its critical points are weak solutions of the Euler--Lagrange equation 
\begin{equation}\label{euler}
\left\{
\begin{array}{ll}
- {\rm div} (a(x,u,\nabla u)) + A_t(x,u,\nabla u) = g(x,u)
  &\hbox{in $\Omega$,}\\
u\ = \ 0 & \hbox{on $\partial\Omega$,}
\end{array}
\right.
\end{equation}
where 
\begin{equation}\label{grad}
\hbox{$A_t(x,t,\xi) = \frac{\partial A}{\partial t}(x,t,\xi),\;
a(x,t,\xi) = (\frac{\partial A}{\partial \xi_1}(x,t,\xi),\dots,\frac{\partial A}{\partial \xi_N}(x,t,\xi))$,}\;
G(x,t) = \int_0^t g(x,\tau) d\tau.
\end{equation}

In addition to the smoothness difficulties which are strictly related to our setting, 
when dealing with unbounded domains the embeddings 
of $X$ in suitable Lebesgue spaces are only continuous,
as it is for the classical problem 
\begin{equation}\label{eulerbase}
- \Delta u + V(x) u = g(x,u),
\qquad u \in W^{1,2}(\R^N),
\end{equation}
whose weak solutions are critical points of $\bar{\J}(u)$
as in \eqref{modello} but with 
$\Omega = \R^N$, $p=2$ and $A_1(x,u) \equiv 1$ and $A_2(x,u) \equiv V(x)$.
Such a lack of compactness
has been overcome in different ways, for example with radial symmetry assumptions
or good hypotheses on potential $V(x)$
(for equation \eqref{eulerbase} see, e.g., \cite{BW, Ra, AS} 
and also \cite{Cerami} and references therein;
for more general quasilinear equations which are model problems of \eqref{euler}
see, e.g., \cite{CS2020, CSS, MS, MS2}).

Here, in order to overcome the lack of
compactness, we make use of an approximating scheme so 
that one positive solution of \eqref{euler} can be found by passing to 
the limit on a sequence $(u_k)_k$ of positive solutions on bounded domains
(a similar approach has been used, for example, in \cite{CDS,LW2,LWW2}).
Furthermore, in stronger assumptions, we prove that either such a solution is nontrivial 
or a diverging sequence of points $(y_k)_k \subset \R^N$ 
exists such that 
\[
\left(\int_{B_1(y_k)} |u_k|^p dx\right)_k \quad \hbox{is far away from zero.}
\]

Anyway, since our main results require a list of technical hypotheses, 
we give their complete statements in Section \ref{main},
while here, in order to highlight at least a model problem and its 
related results, we consider the particular setting of 
the ``modified Schr\"odinger equation'' in 
$\Omega = \R^3$ with $p=2$ and
\[
A(x,t,\xi) \ =\ \frac12 (A_1^*(x) + A_2^*(x) |t|^{s}) |\xi|^2 + \frac12 |t|^2
\quad \hbox{with $A_2^*(x) \not\equiv 0$,} 
\qquad g(x,t)\ =\ |t|^{\mu-2}t,
\]
so that problem \eqref{euler} reduces to
\begin{equation}\label{euler0}
- {\rm div} ((A^*_1(x) + A^*_2(x)|u|^{s}) \nabla u) + 
\frac{s}2 A^*_2(x)\ |u|^{s -2} u\ |\nabla u|^2 + u\ 
=\ |u|^{\mu-2}u \quad\hbox{in $\R^3$.}
\end{equation}

In this setting, the following dichotomy result can be stated.

\begin{theorem}\label{main0}
If a constant $\alpha_0 > 0$ exists such that 
\[
A^*_1,\ A^*_2\ \in\ L^\infty(\R^3)\quad \hbox{and}\quad
A^*_1(x) \ge \alpha_0, \ A^*_2(x) \ge 0 \quad \hbox{a.e. in $\R^3$,}
\]
and if 
\[
3 < 2 + s < \mu< 6,
\]
then a sequence $(u_k)_k$ of bounded positive solutions on bounded domains exists 
which converges a.e. in $\R^3$ to a weak bounded positive solution of problem \eqref{euler0}. \\
Furthermore, either such a solution is nontrivial or 
a constant $\bar\lambda > 0$ and 
a sequence $(y_k)_k \subset \R^N$ exist such that 
\[
|y_k| \to +\infty \qquad \hbox{and}\qquad 
\int_{B_1(y_k)} |u_k|^p dx \ge \bar\lambda \quad \hbox{for all $k \ge 1$.}
\]
\end{theorem}
\medskip

\noindent
\textbf{Notation.$\; $}
Here and in the following, it is $\N = \{1,2,\dots\}$, $|\cdot|$ is the standard norm 
on any Euclidean space as the dimension of the considered vector is clear and no ambiguity arises.

Furthermore, if $C \subseteq \R^N$ is any open domain, we denote by:
\begin{itemize}
\item $B_R(y) = \{x \in \R^N: |x-y| < R\}$ the open ball with center $y$ and radius $R > 0$ in $\R^N$;
\item $B_R = B_R(0)$ for any $R > 0$; 
\item $\meas(C)$ the usual $N$--dimensional Lebesgue measure of $C$;
\item $L^\tau(C)$ the Lebesgue space with norm 
$|u|_{\tau,C} = \left(\int_C|u|^\tau dx\right)^{1/\tau}$ if $1 \le \tau < +\infty$;
\item $L^\infty(C)$ the space of Lebesgue--measurable 
and essentially bounded functions $u :C \to \R$ with norm
\[
|u|_{\infty,C} = \esssup_{C} |u|;
\]
\item $W^{1,p}_0(C)$ and $W^{1,p}(C)$, eventually $W^{1,p}(\R^N)$, the classical Sobolev spaces with
norm 
\[
\|u\|_{W,C} = (|\nabla u|_{p,C}^p + |u|_{p,C}^p)^{\frac1p}\quad \hbox{if $1 \le p < +\infty$.}
\]
\end{itemize}


\section{Variational setting and first properties}
\label{variational}

This section is divided in two parts: firstly, we recall some abstract tools; 
then, we state the variational setting of our problem. 

Let $(X, \|\cdot\|_X)$ be a Banach space with dual space $(X',\|\cdot\|_{X'})$;
furthermore, consider another Banach space $(W,\|\cdot\|_W)$ such that
$X \hookrightarrow W$ continuously,
and a $C^1$ functional $J: X \to \R$.

For simplicity, taking $\beta \in \R$, we say that a sequence
$(u_n)_n\subset X$ is a {\sl Cerami--Palais--Smale sequence at level $\beta$},
briefly {\sl $(CPS)_\beta$--sequence}, if
\[
\lim_{n \to +\infty}J(u_n) = \beta\quad\mbox{and}\quad 
\lim_{n \to +\infty}\|dJ(u_n)\|_{X'} (1 + \|u_n\|_X) = 0.
\]
Moreover, $\beta$ is a {\sl Cerami--Palais--Smale level}, briefly {\sl $(CPS)$--level}, 
if there exists a $(CPS)_\beta$--sequence.

The functional $J$ satisfies the classical Cerami--Palais--Smale condition in $X$ 
at level $\beta$ if every $(CPS)_\beta$--sequence converges in $X$
up to subsequences. Anyway, thinking about the setting of our problem,
in general $(CPS)_\beta$--sequences may also exist which are unbounded in $\|\cdot\|_X$
but converge with respect to $\|\cdot\|_W$ (see \cite[Example 4.3]{CP2017}). 
Then, we can weaken the classical Cerami--Palais--Smale 
condition as follows.  

\begin{definition} 
The functional $J$ satisfies the
{\slshape weak Cerami--Palais--Smale 
condition at level $\beta$} ($\beta \in \R$), 
briefly {\sl $(wCPS)_\beta$ condition}, if for every $(CPS)_\beta$--sequence $(u_n)_n$,
a point $u \in X$ exists such that 
\begin{description}{}{}
\item[{\sl (i)}] $\displaystyle 
\lim_{n \to+\infty} \|u_n - u\|_W = 0\quad$ (up to subsequences),
\item[{\sl (ii)}] $J(u) = \beta$, $\; dJ(u) = 0$.
\end{description}
\end{definition} 

Due to the convergence only in the $W$--norm, the $(wCPS)_\beta$
condition implies that the set of critical points of $J$ at level $\beta$
is compact with respect to $\|\cdot\|_W$; anyway,
this weaker ``compactness'' assumption is enough to prove 
a Deformation Lemma and then some 
abstract critical point theorems (see \cite{CP3}).
In particular, the following generalization of the Mountain Pass Theorem
\cite[Theorem 2.1]{AR} can be stated (for the proof, see \cite[Theorem 1.7]{CP3}
with remarks in \cite[Theorem 2.2]{CS2020}).

\begin{theorem}[Mountain Pass Theorem]
\label{mountainpass}
Let $J\in C^1(X,\R)$ be such that $J(0) = 0$
and the $(wCPS)_\beta$ condition holds for all $\beta \in \R$.
Moreover, assume that two constants
$r$, $\varrho > 0$ and a point $e \in X$ exist such that
\[
u \in X, \; \|u\|_W = r\quad \then\quad J(u) \ge \varrho,
\]
\[
\|e\|_W > r\qquad\hbox{and}\qquad J(e) < \varrho.
\]
Then, $J$ has a critical point $u^* \in X$ such that 
\[
J(u^*) = \inf_{\gamma \in \Gamma} \sup_{\sigma\in [0,1]} J(\gamma(\sigma)) \ge \varrho
\]
with $\Gamma = \{ \gamma \in C([0,1],X):\, \gamma(0) = 0,\; \gamma(1) = e\}$.
\end{theorem}

Now, we state the variational setting for problem \eqref{euler}.
 
From now on, let $\Omega \subseteq \R^N$, $N\ge 2$, 
be an open connected domain with Lipschitz boundary
such that
\begin{equation}\label{infty}
\meas(\Omega) = +\infty,\quad \hbox{eventually $\Omega = \R^N$.}
\end{equation}
For simplicity, we use the notation $|u|_{\tau}$, $|u|_{\infty}$, $\|u\|_{W}$,
instead of $|u|_{\tau,\Omega}$, $|u|_{\infty,\Omega}$, $\|u\|_{W,\Omega}$,
when dealing with the whole domain $\Omega$.

From the Sobolev Embedding Theorems, for any $\tau \in [p,p^*]$
with $p^* = \frac{pN}{N-p}$ if $N > p$, or any $\tau \in [p,+\infty[$ if $p = N$,
the Sobolev space $W_0^{1,p}(\Omega)$ is continuously imbedded in $L^\tau(\Omega)$, 
i.e., $W_0^{1,p}(\Omega) \hookrightarrow L^\tau(\Omega)$ and 
a constant $\sigma_\tau > 0$ exists, such that 
\begin{equation}\label{Sob1}
|u|_\tau\ \le\ \sigma_\tau \|u\|_W \quad \hbox{for all $u \in W_0^{1,p}(\Omega)$}
\end{equation}
(see, e.g., \cite[Corollaries 9.10 and 9.11]{Br}). Clearly, it is $\sigma_p = 1$.
On the other hand, if $p > N$ then  $W_0^{1,p}(\Omega)$ is continuously embedded in $L^\infty(\Omega)$
(see, e.g., \cite[Theorem 9.12]{Br}).\\
In particular, if $p<N$ and $C$ is an open bounded domain,
then a constant $\sigma_* > 0$ exists, independent of $C$ and 
depending only on $p$ and $N$, such that
\begin{equation}\label{stima*}
|v|_{p^*,C}\ \le\ \sigma_* \|v\|_{W,C} \quad \hbox{for all $v \in W_0^{1,p}(C)$.}
\end{equation}
Furthermore, when necessary, if $N=p$
we still consider estimate \eqref{stima*} 
but for any $p^* \in\ ]N, +\infty[$
(see, e.g., \cite[Section 9.4]{Br}). 

Now, for our setting we consider the Banach spaces  
\begin{equation}\label{space}
X := W^{1,p}_0(\Omega) \cap L^\infty(\Omega),\qquad
\|u\|_X = \|u\|_W + |u|_\infty,
\end{equation}
and $W = W^{1,p}_0(\Omega)$ with its norm $\|\cdot\|_W$. 

Clearly, it is 
$X \hookrightarrow W$ continuously and 
we assume $p \le N$ as, otherwise, it is $X = W$ and
the proofs can be simplified.

Reasoning as in \cite[Lemmas 3.1 and 3.2]{CS2020} we have the following results.
 
\begin{lemma}\label{immergo}
For any $p\le \tau \le +\infty$ the Banach space $X$ is continuously imbedded in $L^\tau(\Omega)$, i.e.,
a constant $\sigma_\tau > 0$ exists such that 
\[
|u|_\tau\ \le\ \sigma_\tau \|u\|_X \quad \hbox{for all $u \in X$.}
\]
\end{lemma}

\begin{lemma}\label{immergo2}
If $(u_n)_n \subset X$, $u \in X$, $M > 0$ are such that
\[
|u_n - u|_{p} \to 0 \ \quad\hbox{if $n \to+\infty$,}\qquad
|u_n|_{\infty} \le M\quad \hbox{for all $n \in \N$,}
\]
then $u_n \to u$ also in $L^\tau(\Omega)$ for any $p \le \tau < +\infty$. 
\end{lemma}

For brevity, we introduce the following definition.

\begin{definition}
A map $h : \Omega\times \R^m \to \R$, $m \in \N$, is a $C^k$--Carath\'eodory function, $k \in \N\cup \{0\}$,
if
\[
\begin{split}
&h(\cdot,y) : x \in \Omega\ \mapsto\ h(x,y) \in \R \quad \hbox{is measurable for all $y \in \R^m$,}\\
&h(x,\cdot) : y \in \R^m\ \mapsto\ h(x,y) \in \R \quad \hbox{is $C^k$ for a.e. $x \in \Omega$.}
\end{split}
\]
\end{definition}

Here and in the following, consider $\, A : \Omega \times \R \times \R^N \to \R\,$
and $\, g :\Omega \times \R \to \R\,$ such that, using the notation in \eqref{grad},
the following conditions hold:
\begin{itemize}
\item[$(H_0)$] $A(x,t,\xi)$ is a $C^1$--Carath\'eodory function;
\item[$(H_1)$] a real number $p > 1$ and 
some positive continuous 
functions $\Phi_i$, $\phi_i : \R \to \R$, $i \in \{0,1,2\}$, exist such that
\[
\begin{array}{ccll}
|A(x,t,\xi)| &\le& \Phi_0(t) |t|^p + \phi_0(t)\ |\xi|^p 
& \hbox{a.e. in $\Omega$, for all $(t,\xi) \in \R\times \R^N$,}\\
|A_t(x,t,\xi)| &\le& \Phi_1(t) |t|^{p-1} + \phi_1(t)\ |\xi|^p 
& \hbox{a.e. in $\Omega$, for all $(t,\xi) \in \R\times \R^N$,}\\
|a(x,t,\xi)| &\le& \Phi_2(t) |t|^{p-1} + \phi_2(t)\ |\xi|^{p-1}
& \hbox{a.e. in $\Omega$, for all $(t,\xi) \in \R\times \R^N$;}
\end{array}
\]
\item[$(G_0)$] $g(x,t)$ is a $C^0$--Carath\'eodory function;
\item[$(G_1)$] $a_1$, $a_2 > 0$ and $q \ge p$, $p$ as in $(H_1)$, exist such that
\[
|g(x,t)| \le a_1 |t|^{p-1} + a_2 |t|^{q-1} \qquad
\hbox{a.e. in $\Omega$, for all $t \in \R$.}
\]
\end{itemize}

\begin{remark}\label{remG}
From $(G_0)$--$(G_1)$ it follows that $G(x,t)$ as in \eqref{grad} is
a well defined $C^1$--Ca\-ra\-th\'eo\-do\-ry function in 
$\Omega \times \R$ and
\begin{equation}
\label{alto3}
|G(x,t)| \le \frac{a_1}p |t|^{p} + \frac{a_2}q |t|^{q} \qquad
\hbox{a.e. in $\Omega$, for all $t \in \R$.}
\end{equation}
\end{remark}

\begin{remark}\label{rem0}
From $(H_0)$--$(H_1)$ it follows that
\[
A(x,0,0) = A_t(x,0,0) = 0 \quad \hbox{and}\quad a(x,0,0) = 0 \qquad 
\hbox{for a.e. $x \in \Omega$.}
\] 
Moreover, from $(G_1)$ and \eqref{grad} we have that
\[
G(x,0) = g(x,0) = 0 \qquad 
\hbox{for a.e. $x \in \Omega$.}
\] 
\end{remark}

\begin{remark}\label{rem01}
By comparing assumptions $(H_1)$ and $(G_1)$ with the corresponding 
ones stated in \cite{CP2} when $\Omega$ is bounded, we note that
they are stronger but enough to guarantee the integrability of all 
the involved functions if the Lebesgue measure of $\Omega$ is infinite.\\
In fact, from \eqref{space} and $(H_0)$--$(H_1)$ it follows that
$A(\cdot,u(\cdot),\nabla u(\cdot)) \in L^1(\Omega)$
for all $u \in X$. Furthermore, 
even if no upper bound on $q$ is actually required in $(G_1)$, 
from \eqref{alto3} and Lemma \ref{immergo} it is $G(\cdot,u(\cdot)) \in L^1(\Omega)$
for any $u \in X$, too.\\
On the other hand, $(H_0)$--$(H_1)$ and H\"older inequality
imply that
$a(\cdot,u(\cdot),\nabla u(\cdot)) \cdot \nabla v(\cdot) \in L^1(\Omega)$
and $A_t(\cdot,u(\cdot),\nabla u(\cdot)) v(\cdot) \in L^1(\Omega)$
if $u$, $v \in X$. Moreover, from $(G_0)$--$(G_1)$ 
and, again, H\"older inequality and Lemma \ref{immergo} 
it follows that
$g(\cdot,u(\cdot)) v(\cdot) \in L^1(\Omega)$ if $u$, $v\in X$.
\end{remark}

From Remark \ref{rem01} we can consider the functional $\J : X \to \R$ defined as
\begin{equation}
\label{funct}
\J(u)\ =\ \int_\Omega A(x,u,\nabla u) dx - \int_\Omega G(x,u) dx,
\qquad u \in X,
\end{equation}
and its G\^ateaux differential in $u$ along the direction $v$, with $u$, $v\in X$, as
\begin{equation}
\label{diff}
\langle d\J(u),v\rangle\ =\ \int_{\Omega} a(x,u,\nabla u)\cdot \nabla v\ dx\ +\ 
\int_{\Omega} A_t(x,u,\nabla u) v\ dx\ -\ \int_{\Omega} g(x,u)v\ dx.
\end{equation}

Now, we are ready to state the variational structure of problem \eqref{euler}.

\begin{proposition}\label{smooth1}
Assume that $(H_0)$--$(H_1)$, $(G_0)$--$(G_1)$ are satisfied
for some $p > 1$.
If $(u_n)_n \subset X$, $u \in X$, $M> 0$ are such that 
\[
\|u_n - u\|_W \to 0 \ \quad\hbox{if $n \to+\infty$,}
\]
\[
|u_n|_\infty \le M\quad \hbox{for all $n \in \N$,}
\]
\[
 u_n \to u\quad \hbox{a.e. in $\Omega$} \ \quad\hbox{if $n \to+\infty$,}
\]
then
\[
\J(u_n) \to \J(u)\quad \hbox{and}\quad \|d\J(u_n) - d\J(u)\|_{X'} \to 0
\quad\hbox{if $\ n\to+\infty$.}
\]
Hence, $\J$ is a $C^1$ functional on $X$ with Fr\'echet differential
defined as in \eqref{diff}.
\end{proposition}

\begin{proof}
For simplicity, we set $\J = \J_1 - \J_2$, where
\[
\begin{split}
&\J_1: u \in X \ \mapsto\ \J_1(u) = \int_{\Omega} A(x,u,\nabla u) dx \in \R, \\
&\J_2: u \in X \ \mapsto\ \J_2(u) = \int_{\Omega} G(x,u) dx \in \R,
\end{split}
\]
with related G\^ateaux differentials 
\[
\begin{split}
&\langle d\J_1(u),v\rangle\ =\ \int_{\Omega} a(x,u,\nabla u)\cdot \nabla v\ dx\ +\ 
\int_{\Omega} A_t(x,u,\nabla u) v\ dx,\\
&\langle d\J_2(u),v\rangle\ =\  \int_{\Omega} g(x,u)v\ dx.
\end{split}
\]
By reasoning as in \cite[Proposition 3.10]{MS}
we obtain that  
\[
\J_1(u_n) \to \J_1(u)\quad \hbox{and}\quad \|d\J_1(u_n) - d\J_1(u)\|_{X'} \to 0
\quad\hbox{if $\ n\to+\infty$.}
\]
On the other hand, by reasoning as in the second part of the proof of \cite[Proposition 3.6]{CS2020},
from Lemma \ref{immergo2} with $C=\Omega$ we have that
\[
\J_2(u_n) \to \J_2(u)\quad \hbox{and}\quad \|d\J_2(u_n) - d\J_2(u)\|_{X'} \to 0
\quad\hbox{if $\ n\to+\infty$,}
\]
which complete the proof.
\end{proof}

From Proposition \ref{smooth1} it follows that 
looking for weak (bounded) solutions of \eqref{euler} is equivalent to finding
critical points of the $C^1$ functional $\J$, defined as in \eqref{funct},
on the Banach space $X$ introduced in \eqref{space}.


\section{Set of hypotheses and the main results}
\label{main}

Let $\Omega \subset \R^N$ be an open connected domain with Lipschitz boundary
such that \eqref{infty} holds and
take $A(x,t,\xi)$ and $g(x,t)$ which satisfy $(H_0)$, respectively $(G_0)$.
Then, we consider conditions $(H_1)$ and $(G_1)$ for some $1 < p \le q$,
and introduce the following set of further hypotheses:
\begin{itemize}
\item[$(H_2)$] 
a constant $\alpha_0 > 0$ exists such that
\[
A(x,t,\xi) \ge \alpha_0 (|\xi|^p + |t|^p)
\quad \hbox{a.e. in $\Omega$, for all $(t,\xi) \in \R\times\R^N$;}
\]
\item[$(H_3)$]
a constant $\alpha_1 > 0$ exists such that
\[
a(x,t,\xi)\cdot \xi \ge \alpha_1 |\xi|^p
\quad \hbox{a.e. in $\Omega$, for all $(t,\xi) \in \R\times\R^N$;}
\]
\item[$(H_4)$]
a constant $\eta_0 > 0$ exists such that
\[
A(x,t,\xi) \le \eta_0 (a(x,t,\xi)\cdot \xi + |t|^p)
\quad \hbox{a.e. in $\Omega$, for all $(t,\xi) \in \R\times\R^N$;}
\]
\item[$(H_5)$] 
a constant $\alpha_2 > 0$ exists such that
\[
a(x,t,\xi)\cdot \xi + A_t(x,t,\xi) t \ge \alpha_2 a(x,t,\xi)\cdot \xi
\quad \hbox{a.e. in $\Omega$, for all $(t,\xi) \in \R\times\R^N$;}
\]
\item[$(H_6)$]
some constants $\mu > p$ and $\alpha_3 > 0$ exist such that
\[
\mu A(x,t,\xi) - a(x,t,\xi)\cdot \xi - A_t(x,t,\xi) t \ge \alpha_3 A(x,t,\xi)
\quad \hbox{a.e. in $\Omega$, for all $(t,\xi) \in \R\times\R^N$;}
\]
\item[$(H_7)$]
for all $\xi$, $\xi^* \in \R^N$, $\xi \ne \xi^*$, it is
\[
[a(x,t,\xi) - a(x,t,\xi^*)]\cdot [\xi - \xi^*] > 0 
\quad\hbox{a.e. in $\Omega$, for all $t\in \R$;}
\]
\item[$(G_2)$] taking $\alpha_0$ as in $(H_2)$, then
\[
\exists\ \lim_{t \to 0} \frac{g(x,t)}{|t|^{p-2} t}\ =\ \lambda
\quad \hbox{uniformly for a.e. $x \in \Omega$, with }\; \lambda < \alpha_0 p;
\]
\item[$(G_3)$] taking $\mu$ as in $(H_6)$, then
\[
0 < \mu G(x,t) \le g(x,t) t\qquad\hbox{a.e. in $\Omega$, for all $t \in \R\setminus \{0\}$.}
\]
\end{itemize}

\begin{remark}
Condition $(H_1)$ follows from
assumptions $(H_2)$, $(H_4)$ and the weaker growth hypothesis
\begin{itemize}
\item[$(H'_1)$] a real number $p > 1$ and 
some positive continuous 
functions $\Phi_i$, $\phi_i : \R \to \R$, $i \in \{1,2\}$, exist such that
\[
\begin{array}{ccll}
|A_t(x,t,\xi)| &\le& \Phi_1(t) |t|^{p-1} + \phi_1(t)\ |\xi|^p 
& \hbox{a.e. in $\Omega$, for all $(t,\xi) \in \R\times \R^N$,}\\
|a(x,t,\xi)| &\le& \Phi_2(t) |t|^{p-1} + \phi_2(t)\ |\xi|^{p-1}
& \hbox{a.e. in $\Omega$, for all $(t,\xi) \in \R\times \R^N$.}
\end{array}
\]
\end{itemize}
In fact, for a.e. $x \in\Omega$ and any $(t,\xi) \in \R\times\R^N$,
$(H_2)$, $(H_4)$ and $(H_1')$  imply
\[
0 \le A(x,t,\xi) \le \eta_0 (|a(x,t,\xi)||\xi| + |t|^p)
\le \eta_0 (\Phi_2(t) |t|^{p-1}|\xi| + \phi_2(t)\ |\xi|^{p}+ |t|^p), 
\]
hence, from Young inequality it follows that
\[
0 \le A(x,t,\xi) \le \bar\Phi_0(t) |t|^{p} + \bar\phi_0(t)|\xi|^{p}
\]
for suitable positive continuous 
functions $\bar\Phi_0$, $\bar\phi_0 : \R \to \R$.\\
Anyway, for simplicity, we still consider hypothesis $(H_1)$
even when $(H'_1)$ is enough.
\end{remark}

Now, we point out some direct consequences 
of the previous assumptions.

\begin{remark}\label{sualpha}
From $(H_3)$ and $(H_5)$ it follows that 
$\alpha_2 \le 1$. 
\end{remark}

\begin{remark}\label{conbdd}
From $(H_5)$--$(H_6)$ it follows that 
\[
(\mu - \alpha_3) A(x,t,\xi)\ \ge\ \alpha_2\ a(x,t,\xi)\cdot \xi
\quad \hbox{a.e. in $\Omega$, for all $(t,\xi) \in \R\times\R^N$;}
\]
hence, if also $(H_2)$--$(H_3)$ hold, it is 
$\alpha_3 < \mu$. So, 
\begin{equation}\label{old1}
A(x,t,\xi)\ \ge\ \alpha_4\ a(x,t,\xi)\cdot \xi
\quad \hbox{a.e. in $\Omega$, for all $(t,\xi) \in \R\times\R^N$}
\end{equation}
with $\alpha_4 = \frac{\alpha_2}{\mu - \alpha_3} >0$. Moreover,
from \eqref{old1} and $(H_6)$ we have that
\begin{equation}\label{old2}
\mu A(x,t,\xi) - a(x,t,\xi)\cdot \xi - A_t(x,t,\xi) t\
 \ge\ \alpha_3 \alpha_4\ a(x,t,\xi)\cdot \xi
\quad \hbox{a.e. in $\Omega$, for all $(t,\xi) \in \R\times\R^N$.}
\end{equation}
We note that \eqref{old2} implies \cite[hypothesis $(H_5)$]{CP2}.
\end{remark}

\begin{remark}
We point out that \cite{CP2} deals with functional \eqref{funct}
but when $\Omega$ is a bounded domain; hence, there are some differences 
between the set of hypotheses in \cite{CP2}
and that one here above due to the general failure of the Poincar\'e inequality 
in unbounded domains.
In particular, in our setting, we assume that 
the function $A(x,t,\xi)$ has to be controlled 
by the term $|\xi|^p + |t|^p$ while in 
\cite{CP2} it is enough to consider $|\xi|^p$.
\end{remark}

We note that, from Remark \ref{conbdd}, it is $\alpha_3 < \mu$ in $(H_6)$; thus, the following 
growth estimate can be stated.

\begin{lemma}\label{altoA}
From $(H_0)$, $(H_2)$ and $(H_6)$ it follows that 
\begin{equation}\label{ap1}
A(x,\rho t,\rho\xi) \le\ \rho^{\mu - \alpha_3} \ A(x,t,\xi)\qquad 
\hbox{a.e. in $\Omega$, for all $(t,\xi) \in \R\times\R^N$, $\rho \ge 1$.}
\end{equation}
\end{lemma}

\begin{proof}
Taking any $(t,\xi) \in \R\times\R^N$ and a.e. $x \in \Omega$, 
for all $\rho > 0$ conditions $(H_0)$ and $(H_6)$ imply  
\[
\frac{d}{d\rho}A(x,\rho t,\rho\xi) = A_t(x,\rho t,\rho\xi)t + a(x,\rho t,\rho\xi)\cdot\xi
 \le\ (\mu - \alpha_3)\ \frac 1\rho\ A(x,\rho t,\rho\xi).
\]
Then, from $(H_2)$ we obtain
\[
\frac{\frac{d}{d\rho}A(x,\rho t,\rho\xi)}{A(x,\rho t,\rho\xi)}\ 
 \le\ (\mu - \alpha_3)\ \frac1\rho\quad \hbox{for all $\rho \ge 1$}
\]
and \eqref{ap1} follows by integration with respect to $\rho$.
\end{proof}

\begin{remark}
From $(G_0)$ and $(G_2)$ it follows that
\[
\exists\ \lim_{t \to 0} \frac{G(x,t)}{|t|^{p}}\ =\ \frac{\lambda}{p}
\quad \hbox{uniformly for a.e. $x \in \Omega$.}
\]
Moreover, from assumptions $(G_1)$--$(G_2)$ and direct computations it follows that 
for any $\eps > 0$ a constant $a_\eps > 0$ exists such that
\begin{equation}
\label{altogeps}
|g(x,t)|\ \le\ (\lambda+\eps) |t|^{p-1} + a_\eps |t|^{q-1}\qquad\hbox{a.e. in $\Omega$, 
for all $t \in \R$,}
\end{equation}
and then, if $(G_0)$ holds, 
\begin{equation}
\label{altoGeps}
|G(x,t)|\ \le\ \frac{\lambda + \eps}{p}\ |t|^{p} + \frac{a_\eps}{q} |t|^{q}
\qquad\hbox{a.e. in $\Omega$, for all $t \in \R$.}
\end{equation}
\end{remark}

\begin{remark} \label{bassoG}
From conditions $(G_0)$--$(G_1)$ and $(G_3)$, 
Remark \ref{remG} and direct computations imply that for any 
$\eps > 0$ a function $\eta_\eps \in L^\infty(\Omega)$,
$\eta_\eps(x) > 0$ a.e. in $\Omega$, exists such that
\begin{equation}
\label{basso3}
G(x,t)\ \ge\ \eta_\eps(x)\ |t|^\mu\qquad\hbox{a.e. in $\Omega$, if $|t| \ge \eps$.}
\end{equation}
Hence, from \eqref{alto3}, estimate \eqref{basso3} implies 
\[
p < \mu \le q.
\]
\end{remark}

Now, we are able to state our main existence result.

\begin{theorem} \label{mainthm}
Let $\Omega \subset \R^N$ be an open connected domain with Lipschitz boundary
such that $\meas(\Omega) = +\infty$, eventually $\Omega=\R^N$,
and assume that $A(x,t,\xi)$ and $g(x,t)$ satisfy
conditions $(H_0)$--$(H_7)$, $(G_0)$--$(G_3)$ 
for some $1 < p < \mu \le q$. 
Then, if also
\begin{equation}\label{supq}
q < p^*,
\end{equation} 
a sequence $(u_k)_k$ of bounded positive solutions on bounded domains exists 
which converges a.e. in $\Omega$ to a weak bounded positive solution $\bar{u}$
of problem \eqref{euler}. 
\end{theorem}

Clearly, from Remark \ref{rem0} it follows that
the trivial function $u \equiv 0$, with $\J(0) = 0$,
always satisfies problem \eqref{euler}. 
Unluckly, we have not a direct proof that 
the solution $\bar{u}$ in Theorem \ref{mainthm} is nontrivial
but we are able to obtain more information on 
it if the following stronger assumptions hold:
\begin{itemize}
\item[$(H'_5)$] 
a constant $\alpha_2 > 0$ exists such that
\[
a(x,t,\xi)\cdot \xi + A_t(x,t,\xi) t \ge \alpha_2 (a(x,t,\xi)\cdot \xi + |t|^p)
\quad \hbox{a.e. in $\Omega$, for all $(t,\xi) \in \R\times\R^N$;}
\]
\item[$(G'_2)$] $\quad\displaystyle
\exists\ \lim_{t \to 0} \frac{g(x,t)}{|t|^{p-1}}\ =\ 0\quad$ 
uniformly for a.e. $x \in \Omega$.
\end{itemize}

More precisely, by using a concentration-compactness lemma (see \cite{Li})
we can state the following dichotomy result.

\begin{theorem} \label{mainthm2}
Assume that the hypotheses of Theorem \ref{mainthm} but with 
$(H_5)$ and $(G_2)$ replaced by the stronger conditions
$(H'_5)$ and $(G'_2)$, hold.
Then, either the solution $\bar{u}$ in Theorem \ref{mainthm}
is nontrivial or $\bar\lambda > 0$ and 
a sequence $(y_k)_k \subset \R^N$ exist such that 
\begin{equation}\label{dichotomy}
|y_k| \to +\infty \qquad \hbox{and}\qquad
\int_{B_1(y_k)} |u_k|^p dx \ge \bar\lambda \quad 
\hbox{for all $k \ge 1$.}
\end{equation}
\end{theorem}

Bearing in mind that, in particular, we want to investigate the existence of critical points 
of the model functional $\bar{\J}$ in \eqref{modello},
given a function $A_j : \Omega \times \R \to \R$ and a real number $p > 1$, 
we introduce the following conditions:
\begin{itemize}
\item[$(h_0)$] $A_j(x,t)$ is a $C^1$--Carath\'eodory function;
\item[$(h_1)$] some positive continuous 
functions $\Phi_{j,i} : \R \to \R$, $i \in \{0,1\}$, exist such that
\[
|A_j(x,t)| \le \Phi_{j,0}(t) \quad \hbox{and}\quad
\left|\frac{\partial A_j}{\partial t}(x,t)\right| \le \Phi_{j,1}(t)
\qquad \hbox{a.e. in $\Omega$, for all $t \in \R$;}
\]
\item[$(h_2)$]  
a constant $\alpha_{j,0} > 0$ exists such that
\[
A_j(x,t) \ge \alpha_{j,0} \quad \hbox{a.e. in $\Omega$, for all $t\in \R$;}
\]
\item[$(h_3)$] 
a constant $\alpha_{j,1} > 0$ exists such that
\[
p A_j(x,t) + \frac{\partial A_j}{\partial t}(x,t) t \ge \alpha_{j,1} A_j(x,t)
\quad \hbox{a.e. in $\Omega$, for all $t\in \R$;}
\]
\item[$(h_4)$] 
some constants $\mu > p$ and $\alpha_{j,2} > 0$ exist such that
\[
(\mu - p) A_j(x,t) - \frac{\partial A_j}{\partial t}(x,t) t \ge \alpha_{j,2} A_j(x,t)
\quad \hbox{a.e. in $\Omega$, for all $t\in \R$.}
\]
\end{itemize}

\begin{remark}\label{mod01}
Let $p > 1$ and $A_j(x,t) = A^*_1(x)+ A^*_2(x) B(t)$ for a.e. $x \in \Omega$, for all $t \in \R$.
If some constants $\alpha^*_0$, $\alpha^*_1$, $\alpha^*_2 > 0$ and $\mu > p$
exist such that 
\begin{eqnarray}\label{mod0}
&&A^*_1 \in L^\infty(\Omega)\quad \hbox{and}\quad
A^*_1(x) \ge \alpha^*_0 \quad \hbox{a.e. in $\Omega$,}\\
\label{mod00}
&&A^*_2 \in L^\infty(\Omega)\quad \hbox{and}\quad
A^*_2(x) \ge 0 \quad \hbox{a.e. in $\Omega$,}
\end{eqnarray}
\begin{equation}\label{mod0bis}
\begin{split}
&B \in C^1(\R)\quad \hbox{and}\quad
B(t) \ge 0 \quad \hbox{for all $t \in \R$,}\\
&p B(t) + B'(t) t \ge \alpha^*_1 B(t)
\quad \hbox{and}\quad
(\mu - p) B(t) - B'(t) t \ge \alpha^*_2 B(t)
\quad \hbox{for all $t \in \R$,}
\end{split}
\end{equation}
then $(h_0)$--$(h_4)$ are satisfied. 
\end{remark}

\begin{example}\label{part1}
Let us consider the model case
\begin{equation}\label{mod1}
A(x,t,\xi) = \frac1p A_1(x,t) |\xi|^p + \frac1p A_2(x,t) |t|^p  \quad 
\hbox{a.e. in $\Omega$, for all $(t,\xi) \in \R\times\R^N$,}
\end{equation}
with $p > 1$. Assume that $A_1(x,t)$ satisfies conditions $(h_0)$--$(h_4)$
while $A_2(x,t)$ verifies $(h_0)$, $(h_2)$, $(h_4)$ and, moreover,
\begin{itemize}
\item[$(h_1')$]  a constant $\eta_{2,0} > 0$ and a positive continuous 
function $\Phi_{2,1} : \R \to \R$ exist such that
\[
A_2(x,t) \le \eta_{2,0} \quad \hbox{and}\quad
\left|\frac{\partial A_2}{\partial t}(x,t)\right| \le \Phi_{2,1}(t)
\qquad \hbox{a.e. in $\Omega$, for all $t \in \R$;}
\]
\item[$(h_3')$] 
we have that
\[
p A_2(x,t) + \frac{\partial A_2}{\partial t}(x,t) t \ge 0
\quad \hbox{a.e. in $\Omega$, for all $t\in \R$.}
\]
\end{itemize}
Then, $A(x,t,\xi)$ in \eqref{mod1} verifies the hypotheses $(H_0)$--$(H_4)$,
$(H'_5)$, $(H_6)$--$(H_7)$
as for a.e. $x \in \Omega$ and for all $(t,\xi) \in \R\times\R^N$ there results
\[
A_t(x,t,\xi) = \frac1p \frac{\partial A_1}{\partial t}(x,t)|\xi|^p +
\frac1p \frac{\partial A_2}{\partial t}(x,t)|t|^p + A_2(x,t)|t|^{p-2}t,\quad
a(x,t,\xi) = A_1(x,t) |\xi|^{p-2}\xi.
\]
\end{example}

\begin{example}\label{part2}
In Example \ref{part1} we consider the particular case 
\[
A_1(x,t) = A^*_1(x) + A^*_2(x)\ |t|^s \quad \hbox{and} \quad A_2(x,t) \equiv 1
\qquad \hbox{for a.e. $x \in \Omega$, for all $t \in \R$,} 
\]
i.e., we take the function 
\[
A(x,t,\xi) = \frac1p (A^*_1(x) + A^*_2(x)\ |t|^s) |\xi|^p + \frac1p |t|^p  \quad 
\hbox{a.e. in $\Omega$, for all $(t,\xi) \in \R\times\R^N$,}
\]
with $p > 1$. If \eqref{mod0} and \eqref{mod00} hold and we assume $s > 1$, 
from Remark \ref{mod01} it follows that,
fixing any $\mu > p + s$, $A_1(x,t)$ verifies conditions $(h_0)$--$(h_4)$ 
as $B(t) = |t|^s$ satisfies hypotheses \eqref{mod0bis}. 
\end{example}

\begin{example}\label{part3}
The simplest model function which
satisfies hypotheses $(G_0)$, $(G_1)$, $(G'_2)$, $(G_3)$ is
\[
g(x,t) = |t|^{\mu-2}t \quad \hbox{for a.e. $x \in \R^N$, for all $t \in \R$,
with $q = \mu > p$.}
\]
\end{example}

From Example \ref{part1} we have that Theorem \ref{mainthm2} reduces to the following statement.

\begin{corollary}
Let us consider $p > 1$ and two functions $A_1(x,t)$ and $A_2(x,t)$ such that
$A_1(x,t)$ satisfies conditions $(h_0)$--$(h_4)$
while $A_2(x,t)$ verifies $(h_0)$, $(h_1')$, $(h_2)$--$(h_4)$. 
If $g(x,t)$ verifies hypotheses $(G_0)$, $(G_1)$, $(G'_2)$, $(G_3)$, and 
estimate \eqref{supq} holds,
then for the problem  
\begin{equation}\label{euler1}
\begin{split}
- {\rm div} (A_1(x,u) |\nabla u|^{p-2}\nabla u) &+ 
\frac1p \frac{\partial A_1}{\partial t}(x,u)|\nabla u|^p \\
& +\frac1p \frac{\partial A_2}{\partial t}(x,u)|u|^p + A_2(x,u)|u|^{p-2}u
\ =\ g(x,u)
\qquad\hbox{in $\Omega$,}
\end{split}
\end{equation}
the dicothomy result in Theorem \ref{mainthm2} holds.
\end{corollary}

\begin{remark}
We note that problem \eqref{euler1} is the Euler--Lagrange equation related to the 
model functional $\bar{\J}$ defined in \eqref{modello}.
\end{remark}

More in particular, from Examples \ref{part2} and \ref{part3}
we have that Theorem \ref{mainthm2} reduces to the following statement.

\begin{corollary}
\label{probmod}
Let us consider a function $A^*_1(x)$ verifying \eqref{mod0} and 
a function $A^*_2(x)$ satisfying \eqref{mod00}.
If 
\begin{equation} \label{e0p}
1 < p < p+1 < p + s < \mu < p^*,
\end{equation}
then for the problem  
\begin{equation}\label{euler2}
\begin{split}
- {\rm div} ((A^*_1(x) + A^*_2(x)\ |u|^s) |\nabla u|^{p-2}\nabla u) &+ 
\frac{s}{p} A^*_2(x) |u|^{s - 2} u\ |\nabla u|^p\\
&+ |u|^{p-2}u\ =\ |u|^{\mu-2}u\qquad
\hbox{in $\Omega$,}
\end{split}
\end{equation}
the dicothomy result in Theorem \ref{mainthm2} holds.
\end{corollary}

\begin{remark}
If $1 < p < N$, from \eqref{e0p}
it has to be $N < p^2+p$. In particular, if $p=2$, we obtain 
a solution of \eqref{euler2} when $3 \le N <6$. Thus, 
Corollary \ref{probmod} reduces to Theorem \ref{main0} if $\Omega = \R^3$ and $p=2$.
\end{remark}

As useful in the following, we complete this section pointing out some properties 
of $\J$ in $X$.

\begin{lemma}\label{geo1}
Let us assume that $A(x,t,\xi)$ and $g(x,t)$ 
satisfy hypotheses $(H_0)$--$(H_2)$ and 
$(G_0)$--$(G_2)$ for some 
\begin{equation}\label{subcr}
1 < p < q \le p^*.
\end{equation}
Then, a radius $r > 0$ and a constant $\varrho > 0$ exist such that
\[
u \in X, \quad \|u\|_{W} = r\qquad \then\qquad \J(u) \ge \varrho.
\]
\end{lemma}

\begin{proof}
Taking $\eps > 0$ and any $u \in X$, from $(H_2)$, \eqref{altoGeps} and \eqref{subcr}, 
we have that \eqref{Sob1} implies 
\[
\begin{split}
\J(u) &\ge \alpha_0 \left(\int_\Omega |\nabla u|^p dx 
+ \int_\Omega |u|^p dx\right) - \frac{\lambda+\eps}{p} \int_\Omega |u|^p dx
- \frac{a_\eps}{q} \int_\Omega |u|^q dx\\
&\ge \left(\alpha_0 - \frac{\lambda+\eps}{p}\right) \|u\|_W^p - 
\frac{a_\eps}{q} \sigma_q^q  \|u\|_W^q.
\end{split}
\]
Thus, from $p< q$ and $\lambda < \alpha_0 p$, the thesis follows 
from choosing a constant $\eps$ and a radius $r$ small enough. 
\end{proof}

\begin{lemma}\label{geo2}
Assume that hypotheses $(H_0)$--$(H_2)$, $(H_6)$, 
$(G_0)$--$(G_1)$, $(G_3)$ hold.
Taking a function $\bar{u} \in X$ such that 
\begin{equation}\label{supp}
\meas(C_1) > 0,  
\quad\hbox{with $C_1 = \{x \in \Omega:\ |\bar{u}(x)| \ge 1\}$,}
\end{equation}
then,
\begin{equation}\label{lim1}
\lim_{\rho \to +\infty} \J(\rho \bar{u}) = -\infty.
\end{equation}
\end{lemma}

\begin{proof}
Firstly, we note that from Remark \ref{bassoG} with $\eps=1$, 
a function $\eta_1 > 0$ exists such that   
\[
G(x,\rho \bar{u}(x))\ \ge\ \rho^\mu \eta_1(x)\ |\bar{u}(x)|^\mu
\qquad\hbox{for all $x \in C_1$, for any $\rho \ge 1$.}
\]
Then, from Lemma \ref{altoA} and $(G_3)$ we have that
\[
\begin{split}
\J(\rho \bar{u}) &\le \rho^{\mu -\alpha_3}\ \int_\Omega A(x,\bar{u},\nabla\bar{u}) dx 
- \int_{C_1} G(x,\rho \bar{u}) dx - \int_{\Omega\setminus C_1} G(x,\rho \bar{u}) dx\\
& \le \rho^{\mu -\alpha_3}\ \int_\Omega A(x,\bar{u},\nabla\bar{u}) dx 
- \rho^\mu \int_{C_1} \eta_1(x)\ |\bar{u}|^\mu dx
\end{split}
\]
which implies \eqref{lim1}, as \eqref{supp} implies $\int_{C_1} \eta_1(x)\ |\bar{u}|^\mu dx > 0$. 
\end{proof}

\begin{remark}\label{geo9}
If $R_0 > 0$ is a radius such that 
\begin{equation}\label{erre0}
\meas(\Omega \cap B_{R_0}) > 0,
\end{equation}
then 
a function $\bar{u} \in X$ verifying 
$\supp\ \bar{u} \subset\subset \Omega \cap B_{R_0}$ and
\eqref{supp}, always exists.
In fact, if it enough taking
$\bar{u} = L \varphi_1$ with $L > 0$ large enough and $\varphi_1 > 0$, with
$|\varphi_1|_p = 1$, first eigenfunction of the operator 
$- \Delta_p$ in $W^{1,p}_0(\Omega \cap B_{R_0})$.
We recall that $\varphi_1 \in L^\infty(\Omega \cap B_{R_0})$ (see, e.g.,  \cite{Lin}).
\end{remark}

\begin{remark}\label{geo3}
From Lemma \ref{geo2} and Remark \ref{geo9} 
an element $e \in X$ exists such that
\[
\supp\ e \subset\subset \Omega \cap B_{R_0}, \quad \|e\|_W > r 
\quad \hbox{and}\quad \J(e) < 0,
\]
with $r > 0$ as in Lemma \ref{geo1} and $R_0 >0$
such that \eqref{erre0} is satisfied. \\
Moreover, the function $e$ can be choosen 
so that
\[
e(x) > 0 \quad \hbox{for a.e.}\; x \in \Omega \cap B_{R_0}.
\]
\end{remark}


\section{Approximating problems in bounded domains}
\label{mainbounded}

From now on, let $\Omega \subset \R^N$ be an open connected domain 
with Lipschitz boundary verifying condition \eqref{infty}.
Taking $R_0 > 0$ such that \eqref{erre0} holds,
we consider 
\begin{equation}\label{erre0bis}
R \ge R_0\quad \hbox{and}\quad
\Omega_R = \Omega \cap B_R.
\end{equation}
Thus, for all $R \ge R_0$ we have that
$\meas(\Omega_R) > 0$ and we define 
\begin{equation}\label{icsR}
 X_R = W^{1,p}_0(\Omega_R) \cap L^\infty(\Omega_R) \quad \hbox{equipped with the norm 
$\|u\|_{X_R} = \|u\|_{W,\Omega_R} + |u|_{\infty,\Omega_R}$.}
\end{equation}
Since any function $u \in  X_R$ can be trivially extended 
to $\Omega$ assuming $u(x) = 0$ for all $x \in \Omega \setminus B_R$, then
\[
\|u\|_{W} = \|u\|_{W,\Omega_R}, \qquad  |u|_{\infty} = |u|_{\infty,\Omega_R};
\]
hence, $X_R \hookrightarrow X$ and from Remark \ref{rem0} functional $\J$ in \eqref{funct} is such that
\begin{equation}\label{functR}
\J_R(u) = \J|_{X_R}(u) =
\int_{\Omega_R}A(x,u,\nabla u) dx - \int_{\Omega_R} G(x,u) dx, \quad u \in X_R.
\end{equation}
From \cite[Proposition 3.1]{CP2} we have that $\J_R : X_R \to \R$
is a $C^1$ functional such that
\[
\langle d\J_R(u),v\rangle = \int_{\Omega_R} a(x,u,\nabla u)\cdot \nabla v dx 
+ \int_{\Omega_R} A_t(x,u,\nabla u) v dx- \int_{\Omega_R} g(x,u) v dx, \quad u, v \in X_R.
\]

Now, our aim is to apply Theorem \ref{mountainpass} to $\J_R$ in $X_R$.
Thus, we need to prove the weak Cerami--Palais--Smale condition.

\begin{lemma}\label{cond}
Under the assumptions in Theorem \ref{mainthm},
functional $\J_R$ satisfies the $(wCPS)_\beta$ condition 
for all $\beta \in \R$.
\end{lemma}

\begin{proof}
The result is essentially in \cite[Proposition 4.6]{CP2}, but,
since our assumptions are slightly different,
for the sake of completeness we give here a sketch of the proof. \\
Taking $\beta \in \R$, let $(u_n)_n \subset X_R$ be a $(CPS)_\beta$--sequence, i.e.,
\begin{equation}\label{c1}
\J_R(u_n) \to \beta \quad \hbox{and}\quad \|d\J_R(u_n)\|_{X'_R}(1 + \|u_n\|_{X_R}) \to 0\qquad
\mbox{if $\ n\to+\infty$.}
\end{equation}
We want to prove that $u \in X_R$ exists such that 
\begin{itemize}
\item[$(i)$] $\ \|u_n - u\|_{W,\Omega_R} \to 0$ (up to subsequences), 
\item[$(ii)$] $\ \J_R(u) = \beta$, $d\J_R(u) = 0$.
\end{itemize}
The proof will be divided in the following steps:
\begin{itemize}
\item[1.] $(u_n)_n$ is bounded in $W^{1,p}_0(\Omega_R)$; hence, up to subsequences,
there exists $u \in W^{1,p}_0(\Omega_R)$ such that if $n\to+\infty$, then
\begin{eqnarray}
&&u_n \rightharpoonup u\ \hbox{weakly in $W^{1,p}_0(\Omega_R)$,}
\label{c2}\\
&&u_n \to u\ \hbox{in $L^\tau(\Omega_R)$ for each $\tau \in [1,p^*[$,}
\label{c3}\\
&&u_n \to u\ \hbox{a.e. in $\Omega_R$;}
\label{c4}
\end{eqnarray}
\item[2.] $u \in L^\infty(\Omega_R)$;
\item[3.] having defined $T_m : \R \to \R$ such that
\[
T_m t = \left\{\begin{array}{ll}
t&\hbox{if $|t| \le m$}\\
m \frac t{|t|}&\hbox{if $|t| > m$}
\end{array}\right. ,
\]
if $m \ge |u|_{\infty,\Omega_R} + 1$ then
\begin{equation}\label{p2}
\|d\J_R(T_m u_n)\|_{X'_R} \to 0\quad \hbox{and}\quad
\J_R(T_m u_n) \to \beta \qquad \hbox{as $n \to +\infty$;}
\end{equation}
\item[4.] $\|T_m u_n - u\|_{W,\Omega_R} \to 0$ if $n\to+\infty$; hence, $(i)$ holds; 
\item[5.] $(ii)$ is satisfied.
\end{itemize}
For simplicity, since $\Omega_R$ is an open bounded domain, 
throughout this proof we denote by $\|u\|_{W,\Omega_R}$ the equivalent norm
$|\nabla u|_{p,\Omega_R}$ so that the norm on $X_R$ is replaced by 
the equivalent norm $|\nabla u|_{p,\Omega_R} + |u|_{\infty,\Omega_R}$.
Then, if we compare our set of hypotheses, which imply also \eqref{old2},
with that one in \cite{CP2}, we have that only condition (4.1) in \cite{CP2} does not hold 
but is replaced by $(H_4)$. So, even if many parts of this proof can be obtained 
as in \cite[Proposition 4.6]{CP2}, here, for completeness, we point out when some 
differences occur.
Furthermore, we use the notation $(\eps_n)_n$
for any infinitesimal sequence depending only on $(u_n)_n$ while $(\eps_{m,n})_n$
for any infinitesimal sequence depending not only on $(u_n)_n$ but also on some fixed
integer $m$. Moreover, $b_i$ denotes any strictly positive constant independent of $n$.\\
\emph{Step 1.} From \eqref{functR} -- \eqref{c1}
it follows that
\begin{eqnarray*}
\beta + \eps_n &=& \J_R(u_n) = \int_{\Omega_R} A(x,u_n,\nabla u_n) dx - \int_{\Omega_R} G(x,u_n) dx\\
\eps_n &=& \langle d\J_R(u_n),u_n\rangle\\
 &=& \int_{\Omega_R} a(x,u_n,\nabla u_n)\cdot \nabla u_n dx + \int_{\Omega_R} A_t(x,u_n,\nabla u_n)u_n dx 
 - \int_{\Omega_R} g(x,u_n)u_n dx.
\end{eqnarray*}
Hence, from $(H_2)$, $(H_6)$ and $(G_3)$ we have that
\[
\mu\beta + \eps_n = \mu\J_R(u_n) - \langle d\J_R(u_n),u_n\rangle
\ge \alpha_3 \int_{\Omega_R} A(x,u_n,\nabla u_n) dx \ge \alpha_0\alpha_3 \|u_n\|_{W,\Omega_R}^p,
\]
which implies the boundedness of $(u_n)_n$ in $W^{1,p}_0(\Omega_R)$.
Thus, from the reflexivity of $W^{1,p}_0(\Omega_R)$  
the weak limit \eqref{c2} holds up to subsequences, then \eqref{c3}
and \eqref{c4} follow from Rellich--Kondrachov Theorem (see \cite[Theorem 6.2]{Ad}).\\ 
\emph{Step 2.} The proof is obtained by reasoning as in \cite[Proposition 4.6, \emph{Step 2}]{CP2}.\\
\emph{Step 3.} Define $R_m : \R \to \R$ such that
\begin{equation}\label{resto}
R_mt = t - T_mt = \left\{\begin{array}{ll}
0&\hbox{if $|t| \le m$}\\
\left(1 - \frac{m}{|t|}\right) t &\hbox{if $|t| > m$}
\end{array}\right., \qquad \Omega^{n}_{R,m}= \{x \in \Omega_R: |u_n(x)| > m\}.
\end{equation}
Taking any $m \ge |u|_{\infty,\Omega_R} +1$, we have that 
\[
T_m u = u\quad\hbox{and}\quad R_m u = 0\qquad \hbox{a.e. in $\Omega_R$;}
\]
hence, from \eqref{c3} it follows that as $n \to +\infty$ we obtain
\begin{equation}\label{s11}
T_m u_n \to u \quad\hbox{and}\quad R_m u_n \to 0 \qquad 
\hbox{in $L^\tau(\Omega_R)$ for any $1 \le \tau < p^*$;}
\end{equation}
in particular, $R_m u_n \to 0$ in $L^p(\Omega_R)$ and then
\begin{equation}\label{em11}
\int_{\Omega^n_{R,m}}||u_n| - m|^p dx \to 0,
\end{equation}
while from \eqref{c4} we have that
\begin{equation}\label{pt}
T_mu_{n} \to u\quad\hbox{and}\quad R_m u_n \to 0 
\qquad \hbox{a.e. in $\Omega_R$}
\end{equation}
and also
\begin{equation}\label{bbb8}
\meas(\Omega^n_{R,m}) \to 0.
\end{equation}
Furthermore, by definition, for any $n \in \N$ we have that
\begin{equation}\label{limito}
|T_m u_n|_{\infty,\Omega_R}\ \le\ m 
\end{equation}
and also $\|T_m u_n\|_{X_R} \le \|u_n\|_{X_R}$, 
$\|R_m u_n\|_{X_R} \le \|u_n\|_{X_R}$; 
thus, \eqref{c1} implies
\[
\langle d\J_R(u_n),T_m u_n\rangle \to 0,\qquad \langle d\J_R(u_n),R_m u_n\rangle \to 0
\quad \mbox{as $n \to +\infty$.}
\]
Then, by reasoning as in the proofs of the limits 
(4.17), (4.38) and (4.40) in \cite[Proposition 4.6]{CP2},
we have that
\begin{equation}\label{em1}
\int_{\Omega^n_{R,m}}|\nabla u_n|^p dx \to 0\quad \hbox{as $n \to +\infty$,}
\end{equation}
\begin{equation}\label{emm2}
\int_{\Omega^n_{R,m}} a(x,u_n,\nabla u_n)\cdot \nabla u_n dx \to 0\quad \hbox{as $n \to +\infty$,}
\end{equation}
and also that the first limit in \eqref{p2} is satisfied.\\
Now, in order to prove that the second limit in \eqref{p2} holds,
we note that \eqref{functR} and direct computations imply
\begin{eqnarray*}
\J_R(T_m u_n) &=&
\int_{\Omega_R\setminus{\Omega^n_{R,m}}} A(x,u_n,\nabla u_n)dx +
\int_{\Omega^n_{R,m}} A(x,T_m u_n,0)dx\\
&&-\int_{\Omega_R} G(x,T_m u_n)dx\\
&=& \J_R(u_n) -
\int_{\Omega^n_{R,m}} A(x,u_n,\nabla u_n)dx +
\int_{\Omega^n_{R,m}} A(x,T_m u_n,0)dx\\
&&+\int_{\Omega_R} (G(x,u_n) - G(x,T_mu_n))dx,
\end{eqnarray*}
where $(H_1)$, \eqref{limito} together with \eqref{bbb8} imply
\[
\int_{\Omega^n_{R,m}} A(x,T_mu_n,0)dx \to 0,
\]
while \eqref{c3}, \eqref{s11} and the subcritical growth assumption 
\eqref{alto3} give
\[
\int_{\Omega_R} (G(x,u_n) - G(x,T_m u_n))dx \to 0.
\]
Hence, from \eqref{c1} we have just to prove that 
\begin{equation}\label{eqq1}
\int_{\Omega^n_{R,m}} A(x,u_n,\nabla u_n)dx \to 0.
\end{equation}
To this aim, we note that $(H_2)$ and $(H_4)$ imply that
\[
0 \le \int_{\Omega^n_{R,m}} A(x,u_n,\nabla u_n)dx \le \eta_0 
\int_{\Omega^n_{R,m}} a(x,u_n,\nabla u_n)\cdot \nabla u_n dx 
+ \eta_0 \int_{\Omega^n_{R,m}} |u_n|^p dx,
\]
where by direct computations it is
\[
\int_{\Omega^n_{R,m}} |u_n|^p dx \le 2^{p-1}
\int_{\Omega^n_{R,m}} ||u_n| - m|^p dx + 2^{p-1} m^p \meas(\Omega^n_{R,m}); 
\]
thus, \eqref{eqq1} follows from \eqref{em11}, \eqref{bbb8} and \eqref{emm2}.\\
\emph{Step 4.} 
By using some estimates proved in \emph{Step 3} and hypotheses 
$(H_1)$, $(H_3)$, $(H_7)$, $(G_1)$, 
from an idea introduced in \cite{BMP}
and reasoning as in the proof of \cite[Proposition 4.6, \emph{Step 4}]{CP2},
we prove that 
\[
\int_{\Omega_R}|\nabla T_m u_n - \nabla u|^p dx  \to 0\quad \hbox{if $n\to+\infty$;}
\] 
hence, from \eqref{em1} limit $(i)$ holds.\\
\emph{Step 5.} 
From \eqref{pt}, \eqref{limito} and $(i)$ we have that
 Proposition \ref{smooth1} applies to the uniformly bounded sequence 
$(T_m u_n)_n$, then $(ii)$ follows from \eqref{p2}.
\end{proof}

\begin{proposition}\label{exR}
Under the assumptions in Theorem \ref{mainthm}, 
functional $\J_R$ has at least a critical point $u_R \in X_R$
such that
\[
\varrho\ \le\ \J_R(u_R)\ \le \ \sup_{\sigma\in [0,1]} \J(\sigma e), 
\]
with $\varrho$ as in Lemma \ref{geo1} 
and $e$ as in Remark \ref{geo3}. 
\end{proposition}

\begin{proof}
From Remark \ref{rem0} it follows that $\J_R(0)=0$.
Moreover, from \eqref{functR}, taking $r$, $\varrho > 0$ as in Lemma \ref{geo1},
for any $u \in X_R$ such that $\|u\|_{W,\Omega_R}= r$, 
Lemma \ref{geo1} implies that $\J_R(u) \ge \varrho$.\\
On the other hand, again from \eqref{functR}, taking $e \in X$ as
in Remark \ref{geo3}, condition \eqref{erre0bis} implies $e \in X_R$
so $\|e\|_{W,\Omega_R} > r$ and $\J_R(e) < 0$. Hence, from 
Lemma \ref{cond}, Theorem \ref{mountainpass} applies to $\J_R$
in $X_R$ and a critical point $u_R \in X_R$ exists such that
\[
\J_R(u_R) = \inf_{\gamma \in \Gamma_R} \sup_{\sigma\in [0,1]} \J_R(\gamma(\sigma)) \ge \varrho
\]
with $\Gamma_R = \{ \gamma \in C([0,1],X_R):\, \gamma(0) = 0,\; \gamma(1) = e\}$.\\
Finally, we note that from \eqref{erre0bis} the segment 
\[
\gamma^*: \sigma\in [0,1] \mapsto \sigma e \in X
\]
is in $\Gamma_R$; thus, also the second inequality holds.
\end{proof}
 
\begin{corollary}\label{coroR}
Under the assumptions in Theorem \ref{mainthm}, 
functional $\J_R$ has at least two nontrivial critical points in $X_R$,
one positive and one negative in $\Omega_R$.
\end{corollary}

\begin{proof}
Taking $u \in X$, we have $u = u_+ - u_-$
with $u_\pm(x) = \max\{0,\pm u(x)\}$. Thus,  
let us define 
\[
g_\pm(x,u(x))\ = \ g(x,u_\pm(x))
\qquad\hbox{and then}\qquad G_\pm(x,u(x))\ = \ G(x,u_\pm(x))\quad
\hbox{for a.e. $x\in \Omega$,}
\]
and consider the associated functionals
$\J^\pm : X \to \R$ defined as
\[
\J^\pm(u) \ =\ \int_\Omega A(x,u,\nabla u) dx - \int_\Omega G_\pm(x,u) dx.
\]
From Remark \ref{rem0} and $(G_3)$, we obtain that
\[
G_\pm(x,u(x))\ \le \ G(x,u(x))\quad
\hbox{for a.e. $x\in \Omega$, for all $u \in X$;}
\]
then,
\begin{equation}\label{same}
\J^\pm(u)\ \ge \ \J(u)\quad
\hbox{for all $u \in X$.}
\end{equation}
We note that the statements of Lemmas \ref{geo1} and \ref{geo2}
still hold for $\J^+(u)$ and $\J^-(u)$ in $X$ 
by doing small suitable changes in their proofs.
In particular, from \eqref{same}, we can take the same constants 
$r$, $\varrho > 0$ which are in Lemma \ref{geo1} for both 
$\J^+(u)$ and $\J^-(u)$ in $X$, and can consider $e$ as in Remark 
\ref{geo3} for $\J^+$ while its opposite $-e$ for $\J^-$.\\
Now, by denoting $\J^\pm_R = \J^\pm|_{X_R}$
and reasoning as in the proof of Lemma \ref{cond}, 
we have that both $\J^+_R$ and $\J^-_R$
satisfy the $(wCPS)_\beta$ condition in $X_R$ for all $\beta \in \R$
and then, as in the proof of Proposition \ref{exR},
two points $u^+_R$, $u^-_R \in X_R$ exist
such that 
\[
d \J^+_R (u^+_R) = 0,\quad \J^+_R (u^+_R) \ge \varrho\qquad \hbox{and}\qquad
d \J^-_R (u^-_R) = 0,\quad \J^-_R (u^-_R) \ge \varrho.
\]
At last, the same arguments used at the end of the proof of \cite[Theorem 2.12]{MM}
imply $u^+_R\ge 0$ and $u^-_R \le 0$ a.e. in $\Omega_R$,
and then $u^+_R\ge 0$ is a nontrivial positive, 
respectively $u^-_R$ is a nontrivial negative, critical point of $\J_R$ in $X_R$.
\end{proof}


\section{Proof of Theorem \ref{mainthm}}
\label{mainproof}

In this section, we assume that all the hypotheses of Theorem \ref{mainthm}
are satisfied and, for simplicity, without loss of generality,
we can take $R_0 =1$ in \eqref{erre0}. 
Then, taking any $k \in \N$, we have that
\begin{equation}\label{no0}
\meas(\Omega_k) > 0, \quad \hbox{with $\Omega_k = \Omega \cap B_k$,}
\end{equation}
and, for the radius $R = k$, we can consider the corresponding space $X_k$ 
as in \eqref{icsR} and the related $C^1$ functional 
$\J_k = \J|_{X_k}$ as in \eqref{functR}.

From now on, since any $u \in X_k$, $k \in \N$, can be trivially 
extended as $u = 0$ a.e. in $\Omega \setminus \Omega_k$, 
for simplicity we still denote such an extension by $u$ and then
$X_k \subset X_{k+1}$ for all $k \in \N$.  

\begin{remark}\label{rem00}
From Remark \ref{rem0} we have that if $u\in X_k$ then
\[
\langle d\J_k(u),v\rangle = \langle d\J(u),v\rangle
\quad \hbox{for all $v \in X_k$.}
\]
\end{remark}

So, for each $k \in \N$, \eqref{no0} implies that Proposition \ref{exR}
applies to $\J_k$ in $X_k$ and $u_k \in X_k$ exists such that
\begin{equation}\label{sotto}
\J_k(u_k) = \inf_{\gamma \in \Gamma_k} \sup_{\sigma\in [0,1]} \J_k(\gamma(\sigma)) \ge \varrho
\end{equation}
with $\Gamma_k = \{ \gamma \in C([0,1],X_k):\, \gamma(0) = 0,\; \gamma(1) = e\}$.

We note that from Remark \ref{geo3} it results $\Gamma_k \subset \Gamma_{k+1}$,
thus, 
\begin{equation}\label{cresce}
\J_{k+1}(u_{k+1}) \le \J_k(u_k) \quad \hbox{for all $k \in \N$.} 
\end{equation}

Moreover, from Corollary \ref{coroR} we can choose 
\begin{equation}\label{piu}
u_k(x) \ge 0\quad \hbox{for a.e. $x\in \Omega$, for all $k \in \N$.}
\end{equation}

Then, summing up, from \eqref{piu} and \eqref{sotto}, \eqref{cresce}, Remark \ref{rem00},
a sequence $(u_k)_k \subset X$ of nontrivial positive functions exists such that for every $k \in \N$
it results:
\begin{itemize}
\item[$(P_1)$] $\ u_k \in X_k$ with $u_k = 0$ a.e. in $\Omega \setminus \Omega_k$,
\item[$(P_2)$] $\ \displaystyle \varrho\ \le\ \J(u_k)\ \le \ \J(u_1)$,
\item[$(P_3)$] $\langle d\J(u_k),v\rangle = 0$ for all $v \in X_k$.
\end{itemize}

We note that from $(P_1)$ and $(P_3)$ it follows
\begin{equation}\label{p13}
 \langle d\J(u_k),u_k\rangle = 0
\qquad \hbox{for all $k \in \N$.}
\end{equation}

Now, our aim is proving that the sequence $(u_k)_k$
is bounded in $X$. 
If $p > N$, then it is enough to verify 
that $(\|u_k\|_{W})_k$ is bounded.
 On the contrary, if $p \le N$, the following generalization
of \cite[Theorem  II.5.1]{LU} allows us to prove 
the uniformly boundedness of such a sequence (for the proof,
see \cite[Lemma 5.6]{CSS}).

\begin{lemma}
\label{tecnico} Let $C$ be an open bounded domain in $\R^N$ and
 consider $p$, $q$ so that  $1 < p \le N$
and $p \le q < p^*$ (if $N = p$ we just require that $p^*$ is any
number but in $]q,+\infty[$)
and take $u \in W_0^{1,p}(C)$. If $a^* >0$ and $m_0\in \N$
exist such that 
\[
\int_{C^+_m}|\nabla u|^p dx \le a^* \left(m^q\ \meas(C^+_m) +
\int_{C^+_m} |u|^q dx\right)\quad\hbox{for all $m \ge m_0$,}
\]
with $C^+_m = \{x \in C:\ u(x) > m\}$, then $\displaystyle \esssup_{C} u$
is bounded from above by a positive constant which can be chosen so
that it depends only on $\meas(C^+_{m_0})$, $N$, $p$, $q$, $a^*$, $m_0$, 
$\|u\|_{W,C}$, or better by a positive constant which can be chosen so 
that it depends only on $N$, $p$, $q$, $a^*$, $m_0$ and $a_0^*$
for any $a_0^*$ such that $\max\{\meas(C^+_{m_0}),\|u\|_{W,C}\} \le a_0^*$.
\end{lemma}

Now, we are ready to prove the following crucial properties of  
sequence $(u_k)_k$.

\begin{lemma}\label{bddseq}
A constant $M_0 > 0$ exists such that
\begin{equation}\label{bddX}
\|u_k\|_X \le M_0 \qquad \hbox{for all $k \in \N$.}
\end{equation}
Hence, for all $\tau \ge p$ a constant $M_\tau > 0$ exists 
such that
\begin{equation}\label{bddtau}
|u_k|_{\tau} \le M_\tau
\qquad \hbox{for all $k \in \N$.}
\end{equation}
\end{lemma}

\begin{proof}
From \eqref{p13} and assumptions
$(H_2)$, $(H_6)$ and $(G_3)$ we have that
\[
\mu\J(u_k)\ = \mu\J(u_k) - \langle d\J(u_k),u_k\rangle
\ge \alpha_3 \int_{\Omega} A(x,u_k,\nabla u_k) dx \ge \alpha_0\alpha_3 \|u_k\|_{W}^p,
\]
which, together with $(P_2)$, implies that 
\begin{equation}\label{bddW}
\|u_k\|_{W} \le M_1
\qquad \hbox{for all $k \in \N$.}
\end{equation}
Now, we want to prove that $(u_k)_k$ is bounded in $L^{\infty}(\Omega)$, too.
To this aim, we note that
for a fixed $k \in \N$ either $|u_k|_\infty \le 1$ or $\ |u_k|_\infty > 1$.\\
If $\ |u_k|_\infty > 1$, then from \eqref{piu} it has to be
\[
\esssup_{\Omega} u_k > 1
\]
and then 
\[
\meas(\Omega^{k}_{1}) > 0\quad \hbox{with}\quad
\Omega^{k}_{1} = \{x \in \Omega:\ u_k(x) > 1\}.
\]
From $(P_1)$ it follows that
\[
\Omega^k_{1} \subset \Omega_k,
\]
then $\Omega^{k}_{1}$ is an open bounded domain
such that  
\[
\meas(\Omega^{k}_{1}) \le \int_{\Omega^{k}_{1}} |u_k|^p dx
\le \int_{\Omega} |u_k|^p dx \le \|u_k\|_{W}^p;
\]
hence, from \eqref{bddW} we have 
\begin{equation}\label{bddmeas}
\meas(\Omega^{k}_{1}) \le M_2
\end{equation}
with $M_2 = M_1^p$ independent of $k \in \N$.\\
In order to apply Lemma \ref{tecnico} to $u_k$ on 
the bounded set $\Omega_{k}$,
taking any $m \in \N$, we consider the real function 
$R_m(t)$ in \eqref{resto}, then from \eqref{piu} we have that
\[
R_m u_k(x) = \left\{\begin{array}{ll}
0&\hbox{if $u_k(x) \le m$}\\
u_k(x)-m&\hbox{if $u_k(x) > m$}
\end{array}\right. \qquad \hbox{for a.e. $x \in \Omega$.}
\]
Then, $R_m u_k \in X_k$,
and from $(P_3)$, $(H_5)$, Remark \ref{sualpha}, $(H_3)$
and direct computations it follows that
\begin{eqnarray*}
0 &=& \langle d\J(u_k),R_m u_k\rangle\\
 &=&
\int_{\Omega^{k}_{m}} \left(1-\frac{m}{u_k}\right)\ \Bigl(a(x,u_k,\nabla u_k)\cdot \nabla u_k 
+ A_t(x,u_k,\nabla u_k) u_k \Bigr)dx\\
&&\quad +\int_{\Omega^{k}_{m}} \frac{m}{u_k}\ a(x,u_k,\nabla u_k)\cdot \nabla u_k dx
- \int_{\Omega} g(x,u_k)R_m u_k dx\\
&\ge&\ \alpha_2\ \int_{\Omega^{k}_{m}} 
a(x,u_k,\nabla u_k)\cdot \nabla u_k dx
- \int_{\Omega} g(x,u_k)R_mu_k dx\\
&\ge& \alpha_1\ \alpha_2\ \int_{\Omega^{k}_{m}} |\nabla u_k|^p dx
- \int_{\Omega} g(x,u_k)R_mu_k dx,
\end{eqnarray*}
which implies
\begin{equation}\label{ttu1}
\alpha_1\ \alpha_2\ \int_{\Omega^{k}_{m}} |\nabla u_k|^p dx
\le \int_{\Omega} g(x,u_k)R_mu_k dx,
\end{equation}
with  $\Omega^{k}_{m} = \{x \in \Omega:\ u_k(x) > m\}$.
Obviously, it is
$\Omega^{k}_{m} \subset \Omega^{k}_{1}$.
Since $(G_3)$ implies $g(x,t) > 0$ if $t>0$ for a.e. $x \in \Omega$, then
$g(x,u_k(x)) m > 0$ for a.e. $x \in \Omega_{m}^{k}$; so,
from $(G_1)$ and the Young inequality with exponent $\frac{q}p$ and its conjugate $\frac{q}{q-p}$
(recall that $q>p$ from Remark \ref{bassoG}), it results that
\[
\begin{split}
\int_{\Omega} g(x,u_k)R_m u_k dx &\le 
\int_{\Omega^{k}_{m}} g(x,u_k)u_k dx
\le \int_{\Omega^{k}_{m}} (a_1 u_k^p + a_2 u_k^q) dx\\
&\le \int_{\Omega^{k}_{m}} (\frac{q-p}q\ a_1^{\frac{q}{q-p}} + \frac{p}{q}\ u_k^q) dx 
+ a_2 \int_{\Omega^{k}_{m}} u_k^q dx\\
&= \frac{q-p}q\ a_1^{\frac{q}{q-p}} \meas(\Omega^{k}_{m})
+ \left(\frac{p}{q} + a_2\right) \int_{\Omega^{k}_{m}} u_k^q dx.
\end{split}
\]
Hence, from \eqref{ttu1} we obtain that
\[
\int_{\Omega^{k}_{m}} |\nabla u_k|^p dx \le  
a^* \left(\meas(\Omega^{k}_{m})
+ \int_{\Omega^{k}_{m}} u_k^q dx\right) \quad \hbox{for all $m \ge 1$,}
\]
with $a^*> 0$ independent of $m$ and $k$.
Then, Lemma \ref{tecnico} with $C=\Omega_{k}$ applies and $M_3 > 1$ exists such that
\[
\esssup_{\Omega_{k}} u_k \le M_3
\]
where, from Lemma \ref{tecnico}, \eqref{bddW} and \eqref{bddmeas}
imply that such a constant $M_3$ can be choosen independent of $k \in \N$. \\
At last, \eqref{bddtau} follows from \eqref{bddX}
and Lemma \ref{immergo}.
\end{proof}

\begin{lemma}\label{liminf}
We have that
\[
\liminf_{k\to +\infty} \int_\Omega |u_k|^p dx \ >\ 0.
\]
\end{lemma}

\begin{proof}
Arguing by contradiction, assume that 
\[
\liminf_{k\to +\infty} \int_\Omega |u_k|^p dx \ =\ 0.
\]
Then, up to subsequences, we have that 
\begin{equation}\label{liminf1}
\int_\Omega |u_k|^p dx \ \to \ 0
\end{equation}
and also 
\begin{equation}\label{liminf2}
u_k \ \to \ 0 \qquad \hbox{a.e. in $\Omega$.}
\end{equation}
Furthermore, \eqref{bddX}, \eqref{liminf1} 
and Lemma \ref{immergo2} imply that
\begin{equation}\label{liminf3}
\int_\Omega |u_k|^q dx \ \to \ 0;
\end{equation}
hence, taking any $\eps > 0$, from $(G_3)$ and \eqref{altogeps}
we obtain that
\[
0\ \le\  \int_\Omega g(x,u_k) u_k\ dx\ \le\ (\lambda + \eps) \int_\Omega |u_k|^p dx
+ a_\eps \int_\Omega |u_k|^q dx
\]
which, together with \eqref{liminf1} and \eqref{liminf3}, implies that
\begin{equation}\label{liminf4}
\int_\Omega g(x,u_k) u_k\ dx\ \to\ 0.
\end{equation}
On the other hand, from \eqref{diff}, \eqref{p13}
and \eqref{liminf4} it results
\[
\int_{\Omega} \Bigl(a(x,u_k,\nabla u_k)\cdot \nabla u_k 
+ A_t(x,u_k,\nabla u_k) u_k \Bigr)dx\
=\ \int_{\Omega} g(x,u_k)u_k dx\ \to\ 0,
\]
while, from $(H_5)$ and $(H_3)$, we have that
\[
\int_{\Omega} \Bigl(a(x,u_k,\nabla u_k)\cdot \nabla u_k 
+ A_t(x,u_k,\nabla u_k) u_k \Bigr)dx\
\ge \ \alpha_1\ \alpha_2\ \int_{\Omega} |\nabla u_k|^p dx.
\]
Thus, 
\[
\int_{\Omega} |\nabla u_k|^p dx\ \to\ 0,
\]
which, together with \eqref{liminf1}, gives
\begin{equation}\label{liminf6}
\|u_k\|_W \ \to \ 0.
\end{equation}
At last, from \eqref{bddX}, \eqref{liminf2} and \eqref{liminf6},
Proposition \ref{smooth1} applies and 
\[
\J(u_k) \ \to\ \J(0) = 0
\]
in contradiction with property $(P_2)$.
\end{proof}

From Lemma \ref{bddseq} we have that $\bar{u} \in W^{1,p}_0(\Omega)$ exists
such that
\begin{equation}\label{weak1}
u_k \rightharpoonup \bar{u} \qquad \hbox{weakly in $W^{1,p}_0(\Omega)$,}
\end{equation}
up to subsequences.

\begin{proposition}\label{limp}
A subsequence $(u_{k_n})_n$ exists such that
\begin{equation}\label{strong}
u_{k_n} \to \bar{u} \quad \hbox{in $L^p(\Omega_R)$ for all $R \ge 1$}        
\end{equation}
and 
\begin{equation}\label{limqo}
u_{k_n} \to \bar{u} \quad \hbox{a.e. in $\Omega$.}
\end{equation}
\end{proposition}

\begin{proof}
In order to prove \eqref{strong}, from \eqref{no0} and \eqref{weak1}
we have that
\[
u_k \rightharpoonup \bar{u} \qquad \hbox{weakly in $W^{1,p}(\Omega_1)$,}
\]
where $\Omega_1 = \Omega \cap B_1$ has Lipschitz boundary and $\meas(\Omega_1) > 0$;
so from Rellich--Kondrachov Theorem
a subsequence $(u_{k^1_n})_n$ exists such that
\[
\lim_{n\to+\infty}u_{k^1_n} = \bar{u} \quad \hbox{in $L^{p}(\Omega_1)$}
\]
and
\[
|u_{k^1_n} - \bar{u}|_{p,\Omega_1}\ < \ \frac{1}{2}\quad 
\hbox{for all $n \ge 1$.}
\]
Since, from \eqref{weak1}, it is also 
\[
u_{k^1_n} \rightharpoonup \bar{u} \quad \hbox{weakly in $W^{1,p}(\Omega_2)$,}
\]
again from Rellich--Kondrachov Theorem a subsequence $(u_{k^2_n})_n$ of $(u_{k^1_n})_n$ exists such that
\[
\lim_{n\to+\infty}u_{k^2_n} = \bar{u} \quad \hbox{in $L^{p}(\Omega_2)$}
\]
and
\[
|u_{k^2_n} - \bar{u}|_{p,\Omega_2}\ < \ \frac{1}{2^2}\quad 
\hbox{for all $n \ge 1$.}
\]
By proceeding in an inductive fashion, for $m \ge 2$ we obtain a subsequence 
$(u_{k^{m}_n})_n$ of the sequence $(u_{k^{m-1}_n})_n$
such that 
\[
\lim_{n\to+\infty}u_{k^{m}_n} = \bar{u} \quad \hbox{in $L^{p}(\Omega_{m})$}
\]
and
\[
|u_{k^{m}_n} - \bar{u}|_{p,\Omega_{m}}\ < \ \frac{1}{2^{m}}\quad 
\hbox{for all $n \ge 1$.}
\]
Thus, picking up the `diagonal sequence' $(u_{k_n})_n$ with $u_{k_n} = u_{k_n^n}$, 
we have that
\begin{equation}\label{stimo1}
|u_{k_n} - \bar{u}|_{p,\Omega_{n}}\ < \ \frac{1}{2^{n}}\quad 
\hbox{for all $n \ge 1$.}
\end{equation}
Now, taking a radius $R \ge 1$, from \eqref{stimo1} it follows that
\begin{equation}\label{stimo2}
\int_{\Omega_R}|u_{k_n} - \bar{u}|^p dx\ <\ \frac{1}{2^{n p}}\quad 
\hbox{for all $n \ge R$}
\end{equation}
which implies \eqref{strong}.\\
In order to prove \eqref{limqo}, 
firstly we note that if $n \ge R$ estimate \eqref{stimo2}
implies that
\[
|u_{k_{n+1}} - u_{k_n}|_{p,\Omega_{R}}\ 
\le \ 
|u_{k_{n+1}} - \bar{u}|_{p,\Omega_{R}}
+ |u_{k_n} - \bar{u}|_{p,\Omega_{R}}\ < \ \frac{1}{2^{n +1}} + \frac{1}{2^{n}},
\] 
in other words,
\begin{equation}\label{stimo3}
|u_{k_{n+1}} - u_{k_n}|_{p,\Omega_{R}}\ < \ \frac{3}{2^{n+1}}
\quad \hbox{for all $n \ge R$.}
\end{equation}
Now, reasoning as in the proof of \cite[Theorem 4.8]{Br}, for any $m \ge 1$ and 
a.e. $x \in \Omega$ we set 
\[
g_m(x) = \sum_{n=1}^m |u_{k_{n+1}}(x) - u_{k_n}(x)|.
\]
Thus, by definition, we have that
\begin{equation}\label{stimo3bis}
0\ \le\ g_m(x) \ \le g_{m+1}(x) \quad\hbox{for all $m \ge 1$} 
\end{equation}
which implies that  
\begin{equation}\label{stimo4}
g(x) = \lim_{m\to +\infty} g_m(x) = \sup_{m} g_m(x) \quad \hbox{with}
\quad 0 \le g(x) \le +\infty
\end{equation}
and also
\begin{equation}\label{stimo4ter}
|g_m|_{p,\Omega_R} \ \le\ |g_{m+1}|_{p,\Omega_R} \quad \hbox{for all $m \ge 1$.}
\end{equation}
We claim that 
\begin{equation}\label{stimo4bis}
g(x) < +\infty \quad \hbox{for a.e. $x \in \Omega$.}
\end{equation}
To this aim, we note that from \eqref{bddX} it results
\[
\begin{split}
|g_m|_p\ &\le\ \sum_{n=1}^m |u_{k_{n+1}} - u_{k_n}|_p\ \le\
\sum_{n=1}^m \left(|u_{k_{n+1}}|_p + |u_{k_n}|_p\right) \\
&\le\ \sum_{n=1}^m \left(\|u_{k_{n+1}}\|_X + \|u_{k_n}\|_X\right)
\ \le\ 2 m M_0 \ <\ +\infty;  
\end{split}
\]
hence, $g_m \in L^p(\Omega)$ and also
\begin{equation}\label{stimo5}
|g_m|_{p,\Omega_R} \ \le\ 2 m M_0 \quad \hbox{for all $R \ge 1$.}
\end{equation}
On the other hand, fixing any $R \ge 2$ and taking $\bar{m} \ge R$, 
for all $m \ge \bar{m}$, direct computations, \eqref{stimo3} and \eqref{stimo5} 
imply that
\[
\begin{split}
|g_m|_{p,\Omega_R}\ &=\ \left|\sum_{n=1}^{\bar{m}-1} |u_{k_{n+1}} - u_{k_n}|
\ + \ \sum_{n=\bar{m}}^m |u_{k_{n+1}} - u_{k_n}|\right|_{p,\Omega_R}\\
&\le\  |g_{\bar m -1}|_{p,\Omega_R}\
+ \ \sum_{n=\bar{m}}^m |u_{k_{n+1}} - u_{k_n}|_{p,\Omega_R}\\
&\le\ 2 (\bar{m} - 1) M_0 + \sum_{n=\bar{m}}^m \frac{3}{2^{n+1}}\
\le\ 2 (\bar{m} - 1) M_0 + \frac{3}{2} \sum_{n=0}^{+\infty} \frac{1}{2^{n}} = 2 (\bar{m} - 1) M_0 + 3. 
\end{split}
\]
Thus,  
\begin{equation}\label{stimo6}
\sup_{m}|g_m|_{p,\Omega_R} < +\infty.
\end{equation}
Now, from \eqref{stimo4} together with \eqref{stimo3bis}, 
Beppo--Levi Theorem and \eqref{stimo4ter}
we obtain that
\[
\int_{\Omega_R} |g|^{p} dx\ =\ \lim_{m\to+\infty} \int_{\Omega_R} |g_m|^{p} dx\
=\ \sup_{m}|g_m|^p_{p,\Omega_R}.
\]
Hence, \eqref{stimo6} implies that
\[
\int_{\Omega_R} |g|^{p} dx\ < +\infty
\]
and then
\[
g(x) < +\infty \quad \hbox{for a.e. $x \in \Omega_R$.}
\]
Thus, \eqref{stimo4bis} follows from the arbitrariness of $R\ge 2$.\\
The next step of this proof is showing that
the sequence $(u_{k_n}(x))_n$ is convergent for a.e. $x \in \Omega$. 
Indeed, taking any $n$, $j \in \N$, we have that
\[
\begin{split}
|u_{k_{n+j}}(x) - u_{k_n}(x)|\ &\le \ \sum_{i=n}^{n+j-1} |u_{k_{i+1}}(x) - u_{k_{i}}(x)|\\
&=\ \sum_{i=1}^{n+j-1} |u_{k_{i+1}}(x) - u_{k_{i}}(x)| 
- \sum_{i=1}^{n-1} |u_{k_{i+1}}(x) - u_{k_{i}}(x)| \\
&=\ g_{n+j-1}(x) - g_{n-1}(x),
\end{split}
\]
where, from \eqref{stimo4} and \eqref{stimo4bis}, 
for a.e. $x \in \Omega$ the sequence $(g_m(x))_m$
is converging and, then, Cauchy. 
Thus, for a.e. $x \in \Omega$ the sequence $(u_{k_n}(x))_n$ is Cauchy and, then, convergent.\\
Finally, taking any $R \ge 1$ 
from \eqref{strong} it has to be
\[
\lim_{n\to +\infty} u_{k_n}(x) = \bar{u}(x) \quad \hbox{for a.e. $x \in \Omega_R$};
\]
so, \eqref{limqo} follows from the arbitrariness of $R \ge 1$.
\end{proof}

The same arguments developed in the proof of Proposition \ref{limp}, 
allow us to state the following result.

\begin{proposition}\label{qo}
Let $(\omega_n)_n \subset L^p(\Omega)$ be a bounded sequence
such that 
\[
\omega_n \to \omega \quad \hbox{in $L^p(\Omega_R)$ for all $R \ge 1$,}
\]
for a certain $\omega \in L^p(\Omega)$. Then, up to a subsequence, 
\[
\omega_n \to \omega \quad \hbox{a.e. in $\Omega$.}
\]
\end{proposition}

From now on, we will consider the subsequence $(u_{k_n})_n$ as in Proposition \ref{limp}
and, for simplicity, we denote it as $(u_k)_k$ so 
\eqref{strong}, respectively \eqref{limqo}, can be written as
\begin{equation}\label{strong1}
u_{k} \to \bar{u} \quad \hbox{in $L^p(\Omega_R)$ for all $R \ge 1$,}
\end{equation}
respectively
\begin{equation}\label{limqo2}
u_{k} \to \bar{u} \quad \hbox{a.e. in $\Omega$.}
\end{equation} 
Clearly, from \eqref{piu} we have that 
\[
\bar{u}(x) \ge 0 \quad \hbox{for a.e. $x \in \Omega$.}
\]

We are ready to prove some further properties
on this sequence $(u_k)_k$ and its limit $\bar{u}$.
 
\begin{lemma}\label{bddseq1}
$\bar{u} \in L^\infty(\Omega)$; hence, $\bar{u} \in X$.
\end{lemma}

\begin{proof}
From \eqref{limqo2} a subset $Z \subset \Omega$ exists
such that $|Z|= 0$ and 
\[
u_k(x) \to \bar{u}(x)
\quad \hbox{for all $x\in \Omega\setminus Z$.}
\] 
Then, for each $x\in \Omega\setminus Z$ an integer $\nu_{x,1}\in\N$ exists such that
\[
|u_k(x) - \bar{u}(x)| < 1\quad\mbox{ for all } k\ge\nu_{x,1},
\]
and, in particular, from \eqref{bddX} we have that
\[
|\bar{u}(x)| \le |\bar{u}(x)- u_{\nu_{x,1}}(x)| +|u_{\nu_{x,1}}(x)| \le 1+ M_0.
\]
Hence, $|\bar{u}(x)| \le 1 + M_0$ for all $x\in \Omega\setminus Z$ 
which completes the proof.
\end{proof}

\begin{remark}\label{limtau}
From Lemma \ref{bddseq}, Lemma \ref{bddseq1} and \eqref{strong1},
taking any $1 \le \tau < +\infty$ 
we have that
\[
u_{k} \to \bar{u} \quad \hbox{in $L^\tau(\Omega_R)$, for all $R \ge 1$.}
\]
\end{remark}

\begin{proposition}\label{limwpr}
We have that
\begin{equation}\label{strongwpr}
u_k \to \bar{u} \quad \hbox{strongly in $W^{1,p}(\Omega_R)$ for all $R \ge 1$.}
\end{equation}
Furthermore, up to a subsequence, it results
\begin{equation}\label{limgradqo}
\nabla u_{k} \to \nabla \bar{u} \quad \hbox{a.e. in $\Omega$.}
\end{equation}
\end{proposition}

\begin{proof}
The main tools of the proof are essentially as in the proof of \textsl{Step 4}
in \cite[Proposition 4.6]{CP2}, but
since here some dissimilarities occur, we give more details. 
For simplicity, throughout this proof we use the notation $(\eps_k)_k$ 
for any infinitesimal sequence which depends on $(u_k)_k$
and $b_i$ for any strictly positive constant independent of $k$.\\
Taking any $R \ge 1$, from \eqref{strong1} it is enough to prove that
\begin{equation}\label{strong2}
|\nabla u_{k} - \nabla \bar{u}|_{p,\Omega_R} \to 0.
\end{equation}
To this aim, defining $v_{k}=u_k - \bar{u}$, we have that 
\eqref{weak1} and \eqref{limqo2} imply 
\begin{equation}\label{cc2}
v_{k} \rightharpoonup 0\quad \hbox{weakly in $W^{1,p}_0(\Omega)$,}
\end{equation}
respectively
\begin{equation}\label{cc22}
v_{k} \to 0 \quad\hbox{a.e. in $\Omega$,}
\end{equation}
while from \eqref{strong1} it follows that
\[
v_{k} \to 0 \quad \hbox{in $L^p(\Omega_{R+1})$.} 
\]
Moreover, from \eqref{bddX} and Lemma \ref{bddseq1} it is
\begin{equation}\label{cc23}
v_k \in X\quad \hbox{and}\quad |v_{k}|_\infty \le \bar{M}_0 \quad\hbox{for all $k \in\N$}
\end{equation}
with $\bar{M}_0 = M_0 + |\bar{u}|_\infty$.\\
Now, following an idea introduced in \cite{BMP}, 
let us consider the real map 
\[
\psi(t) = t \e^{\eta t^2},
\] 
where $\eta > (\frac{\beta_2}{2\beta_1})^2$ will be fixed once
$\beta_1$, $\beta_2 > 0$ are chosen in a suitable way later on. 
By definition,
\begin{equation}\label{eq4}
\beta_1 \psi'(t) - \beta_2 |\psi(t)| > \frac{\beta_1} 2\qquad \hbox{for all $t \in \R$.}
\end{equation}
We note that \eqref{cc23} gives
\begin{equation}\label{stim1}
|\psi(v_{k})| \le \psi(\bar{M}_0),\quad 0<\psi'(v_{k}) \le \psi'(\bar{M}_0) \qquad\hbox{a.e. in $\Omega$,}
\end{equation}
while \eqref{cc22} implies that
\begin{equation}\label{stim2}
\psi(v_{k}) \to 0, \quad
\psi'(v_{k}) \to 1 \qquad\hbox{a.e. in $\Omega$.}
\end{equation}
Furthermore, let $\chi_R \in C^{\infty}(\Omega)$ be a cut--off function 
such that
\begin{equation}\label{cut1}
\chi_R(x) = \left\{\begin{array}{ll}
1 &\hbox{if $|x| \le R$}\\
0 &\hbox{if $|x| \ge R + 1$}
\end{array}\right.,\qquad
\hbox{with}\quad 0 \le \chi_R(x) \le 1\quad \hbox{for all $x \in \Omega$,}
\end{equation}
and
\begin{equation}\label{cut2}
|\nabla \chi_R(x)| \le 2\quad \hbox{for all $x \in \Omega$.}
\end{equation}
Thus, for every $k \in \N$ we consider the new function
\[
w_{R,k}: x \in \Omega \mapsto w_{R,k}(x) = \chi_R(x) \psi(v_k(x)) \in \R.
\]
From \eqref{cut1}, together with \eqref{stim1}, we have that
$\supp w_{R,k} \subset \supp \chi_R \subset \Omega_{R+1}$ and  
\begin{equation}\label{stim3}
|w_{R,k}(x)| \le \psi(\bar{M}_0) \qquad\hbox{a.e. in $\Omega$;}
\end{equation}
moreover, \eqref{stim2} implies that
\begin{equation}\label{stim4}
w_{R,k} \to 0 \qquad\hbox{a.e. in $\Omega$,}
\end{equation}
while from
\begin{equation}\label{grad1}
\nabla w_{R,k} = \psi(v_{k}) \nabla \chi_R + \chi_R \psi'(v_{k}) \nabla v_k\qquad\hbox{a.e. in $\Omega$}
\end{equation}
and \eqref{cc23}, \eqref{stim1}, \eqref{cut1}, \eqref{cut2} 
it follows that $w_{R,k} \in X_{R+1}$.
Hence, for all $k \ge R+1$ we have that
\[
w_{R,k} \in X_{k}
\]
so $(P_3)$, \eqref{diff} and \eqref{grad1} imply that
\begin{equation}\label{zero}
\begin{split}
0\ =\ &\langle d\J(u_k),w_{R,k}\rangle\\
 =\ & \int_{\Omega_{R+1}} \psi(v_{k})\ a(x,u_k,\nabla u_k) \cdot \nabla \chi_R dx
+ \int_{\Omega_{R+1}} \chi_R\ \psi'(v_{k})\ a(x,u_k,\nabla u_k) \cdot \nabla v_{k} dx\\
&+\int_{\Omega_{R+1}}A_t(x,u_k,\nabla u_k)w_{R,k} dx -\int_{\Omega_{R+1}} g(x,u_k) w_{R,k} dx.
\end{split}
\end{equation}
We note that $(G_0)$--$(G_1)$ together with \eqref{bddX}, \eqref{limqo2}, \eqref{stim3}, 
\eqref{stim4} and Lebesgue dominated convergence theorem on the 
bounded set $\Omega_{R+1}$ imply that
\begin{equation}\label{ip2}
\int_{\Omega_{R+1}} g(x,u_k) w_{R,k} dx \to 0.
\end{equation}
Moreover, from $(H_1)$, \eqref{bddX}, \eqref{cut2}, H\"older inequality and direct computations 
it follows that
\[
\begin{split}
\int_{\Omega_{R+1}} \big|\psi(v_{k})\ a(x,u_k,\nabla u_k) \cdot \nabla \chi_R\big| dx &\le
b_1 \int_{\Omega_{R+1}} |\psi(v_{k})| dx + b_2 \int_{\Omega_{R+1}} |\nabla u_k|^{p-1} |\psi(v_{k})| dx\\
&\le 
b_1 \int_{\Omega_{R+1}} |\psi(v_{k})| dx + b_2 \|u_k\|_W^{p-1} \left(\int_{\Omega_{R+1}} |\psi(v_{k})|^p dx\right)^{\frac1p},
\end{split}
\]
where \eqref{stim1}, \eqref{stim2} and Lebesgue dominated convergence theorem on the 
bounded set $\Omega_{R+1}$ give
\begin{equation}\label{uno0}
\int_{\Omega_{R+1}} |\psi(v_{k})| dx \to 0\quad
\hbox{and}\quad \int_{\Omega_{R+1}} |\psi(v_{k})|^p dx \to 0,
\end{equation}
so from \eqref{bddX} we have that
\begin{equation}\label{uno}
\int_{\Omega_{R+1}} \psi(v_{k})\ a(x,u_k,\nabla u_k) \cdot \nabla \chi_R dx \to 0.
\end{equation}
On the other hand, $(H_1)$, \eqref{bddX}, \eqref{cut1}, 
\eqref{uno0}, $(H_3)$ and direct computations imply that
\begin{equation}\label{due}
\begin{split}
\int_{\Omega_{R+1}}\left|A_t(x,u_k,\nabla u_k)w_{R,k}\right| dx
&\le b_3 \int_{\Omega_{R+1}} |\psi(v_{k})| dx + b_4 \int_{\Omega_{R+1}} \chi_R |\nabla u_k|^p|\psi(v_{k})| dx\\
& \le \eps_k + b_5\ \int_{\Omega_{R+1}} \chi_R\ |\psi(v_{k})|\ a(x,u_k,\nabla u_k)\cdot \nabla u_k dx.
\end{split}
\end{equation}
We note that
\[
\begin{split}
\int_{\Omega_{R+1}} \chi_R\ |\psi(v_{k})|\ a(x,u_k,\nabla u_k)\cdot \nabla u_k dx
=\ &\int_{\Omega_{R+1}} \chi_R\ |\psi(v_{k})|\ a(x,u_k,\nabla u_k)\cdot \nabla v_k dx\\
&\; + \int_{\Omega_{R+1}} \chi_R\ |\psi(v_{k})|\ a(x,u_k,\nabla u_k)\cdot \nabla \bar{u}\ dx
\end{split}
\]
where $(H_1)$, \eqref{bddX}, \eqref{cut1} and H\"older inequality imply that
\[
\begin{split}
&\int_{\Omega_{R+1}} \chi_R\ |\psi(v_{k})|\ |a(x,u_k,\nabla u_k)\cdot \nabla \bar{u}| dx
\le \int_{\Omega_{R+1}} |\psi(v_{k})| (b_6  + b_7 |\nabla u_k|^{p-1}) |\nabla \bar{u}| dx\\
&\qquad \le b_6 \int_{\Omega_{R+1}} |\psi(v_{k})| |\nabla \bar{u}| dx +
b_7 \|u_k\|_W^{p-1}
\left(\int_{\Omega_{R+1}} |\psi(v_{k})|^p\ |\nabla \bar{u}|^p dx\right)^{\frac1p}
\end{split}
\]
with 
\[
\int_{\Omega_{R+1}} |\psi(v_{k})| |\nabla \bar{u}| dx \to 0\quad\hbox{and}\quad
\int_{\Omega_{R+1}} |\psi(v_{k})|^p\ |\nabla \bar{u}|^p dx \to 0
\]
thanks, again, to \eqref{stim1}, \eqref{stim2} and Lebesgue dominated convergence theorem on the 
bounded set $\Omega_{R+1}$. Thus, from \eqref{bddX} estimate \eqref{due} becomes
\begin{equation}\label{tre}
\int_{\Omega_{R+1}}\left|A_t(x,u_k,\nabla u_k)w_{R,k}\right| dx
 \le \eps_k + b_5\ \int_{\Omega_{R+1}} \chi_R\ |\psi(v_{k})|\ a(x,u_k,\nabla u_k)\cdot \nabla v_k dx.
\end{equation}
Then, by using \eqref{ip2}, \eqref{uno} and \eqref{tre} in \eqref{zero} we obtain 
\begin{equation}\label{quattro}
\eps_{k} \ge \int_{\Omega_{R+1}} \chi_R\ h_{k}\ a(x,u_k,\nabla u_k) \cdot \nabla v_{k} dx,
\end{equation}
with 
\[
h_{k} (x) = \psi'(v_{k}(x)) - b_5 |\psi(v_{k}(x))|\quad \hbox{a.e. in $\Omega$.}
\]
We note that \eqref{stim1} and \eqref{stim2} give
\begin{equation}\label{stim10}
h_{k}(x) \to 1\quad\hbox{a.e. in $\Omega$} \quad 
\hbox{and}\quad
|h_{k}(x)| \le \psi'(\bar{M}_0) + b_{5}|\psi(\bar{M}_0)|\quad\hbox{a.e. in $\Omega$,}
\end{equation}
and also, if in the definition of $\psi$ we fix 
$\beta_1 = 1$ and $\beta_2 = b_5$, \eqref{eq4} implies
\begin{equation}\label{stim5}
h_{k}(x) > \frac12 \quad\hbox{a.e. in $\Omega$, for all $k \in \N$.}
\end{equation}
At last, we have that
\[
\begin{split}
&\int_{\Omega_{R+1}} \chi_R\ h_{k}\ a(x,u_k,\nabla u_k) \cdot \nabla v_{k} dx
=\  \int_{\Omega_{R+1}} \chi_R\ h_{k}\ \big(a(x,u_k,\nabla u_k) - a(x,u_k,\nabla\bar{u})\big) \cdot \nabla v_{k} dx\\
&\qquad + \int_{\Omega_{R+1}} \chi_R\ \big(h_{k}\ a(x,u_k,\nabla \bar{u}) 
- a(x,\bar{u},\nabla\bar{u})\big) \cdot \nabla v_{k} dx
+  \int_{\Omega_{R+1}} \chi_R\ a(x,\bar{u},\nabla\bar{u}) \cdot \nabla v_{k} dx
\end{split}
\]  
where \eqref{cc2} implies
\[
\int_{\Omega_{R+1}} \chi_R\ a(x,\bar{u},\nabla\bar{u}) \cdot \nabla v_{k} dx \to 0.
\]
On the other hand, \eqref{cut1} and H\"older inequality give
\[
\begin{split}
&\int_{\Omega_{R+1}} \chi_R\ \big|\big(h_{k}\ 
a(x,u_k,\nabla \bar{u}) - a(x,\bar{u},\nabla\bar{u})\big) \cdot \nabla v_{k}\big| dx\\
&\qquad\qquad \le\ \|v_{k}\|_W\
\left(\int_{\Omega_{R+1}} \big|h_{k}\ a(x,u_k,\nabla \bar{u}) 
- a(x,\bar{u},\nabla\bar{u})\big|^{\frac{p}{p-1}} dx\right)^{\frac{p-1}{p}},
\end{split}
\]
where $(H_0)$, \eqref{limqo2} and \eqref{stim10} imply
\[
h_{k}\ a(x,u_k,\nabla \bar{u}) - a(x,\bar{u},\nabla\bar{u}) \to 0 \quad\hbox{a.e. in $\Omega$}
\]
with 
\[
|h_{k}\ a(x,u_k,\nabla \bar{u}) - a(x,\bar{u},\nabla\bar{u})|^{\frac{p}{p-1}} \le b_8 + b_9 |\nabla\bar{u}|^p 
\quad\hbox{a.e. in $\Omega$}
\]
thanks to $(H_1)$, \eqref{bddX}, Lemma \ref{bddseq1} and \eqref{stim10},
so Lebesgue dominated convergence theorem on the 
bounded set $\Omega_{R+1}$ applies and, again, from \eqref{bddX}
we obtain that
\[
\int_{\Omega_{R+1}} \chi_R\ \big|\big(h_{k}\ 
a(x,u_k,\nabla \bar{u}) - a(x,\bar{u},\nabla\bar{u})\big) \cdot \nabla v_{k}\big| dx\to 0.
\]
Thus, by using all the previous estimates in \eqref{quattro},
we note that \eqref{cut1}, \eqref{stim5} and $(H_7)$ imply
both 
\[
 \chi_R\ h_{k}\ \big(a(x,u_k,\nabla u_k) - a(x,u_k,\nabla\bar{u})\big) \cdot \nabla v_{k} \ge 0
\quad\hbox{a.e. in $\Omega$}
\]
and
\[
\begin{split}
\eps_{k} &\ge \int_{\Omega_{R+1}} \chi_R\ h_{k}\ \big(a(x,u_k,\nabla u_k) 
- a(x,u_k,\nabla\bar{u})\big) \cdot \nabla v_{k} dx\\
& \ge\ \frac{1}{2}\ \int_{\Omega_{R}} \chi_R\ \big(a(x,u_k,\nabla u_k) 
- a(x,u_k,\nabla\bar{u})\big) \cdot \nabla v_{k} dx\\
& =\ \frac{1}{2}\ \int_{\Omega_{R}} \big(a(x,u_k,\nabla u_k) 
- a(x,u_k,\nabla\bar{u})\big) \cdot \nabla v_{k} dx\ \ge\ 0
\end{split}
\]
which implies
\[
\int_{\Omega_{R}} \big(a(x,u_k,\nabla u_k) 
- a(x,u_k,\nabla\bar{u})\big) \cdot \nabla v_{k} dx\ \to\ 0.  
\]
Finally, hypotheses $(H_1)$, $(H_3)$, $(H_7)$, limit estimates
\eqref{weak1}, \eqref{limqo2} and the uniform boundedness 
in the $L^\infty$--norm \eqref{bddX}
allow us to apply \cite[Lemma 5]{BMP}, thus \eqref{strong2} holds.\\
At last, we note that \eqref{limgradqo} follows from \eqref{strong2}
and Proposition \ref{qo}.
\end{proof}

\begin{proof}[Proof of Theorem \ref{mainthm}]
We claim that $d\J(\bar{u})= 0$ in $X$.  
To this aim, from the density
of $C^\infty_c(\Omega)$ in $X$ it is enough to prove that
\begin{equation}\label{crit0}
\langle d\J(\bar{u}),\varphi\rangle = 0 \quad \hbox{for all $\varphi \in C_c^\infty(\Omega)$}
\end{equation}
with $C_c^\infty(\Omega) = \{\varphi \in C^\infty(\Omega):\ \supp\varphi\subset\subset \Omega\}$.\\
Taking $\varphi \in C_c^\infty(\Omega)$, a radius $R \ge 1$ exists
such that $\supp \varphi \subset \Omega_R$. Thus, for all $k\ge R$
we have that $\varphi \in X_k$ so $(P_3)$ applies and
$\langle d\J(u_k),\varphi\rangle = 0$ implies that
\begin{equation}\label{crit1}
\begin{split}
0\ \le\ |\langle d\J(\bar{u}),\varphi\rangle| &=
|\langle d\J(u_k),\varphi\rangle - \langle d\J(\bar{u}),\varphi\rangle|\\
&\le \int_{\Omega_R}\big|a(x,u_k,\nabla u_k) - a(x,\bar{u},\nabla \bar{u})\big| |\nabla\varphi| dx\\
&\qquad+ |\varphi|_\infty \int_{\Omega_R}\big|A_t(x,u_k,\nabla u_k) - A_t(x,\bar{u},\nabla \bar{u})\big| dx\\
&\qquad+ |\varphi|_\infty \int_{\Omega_R}\big|g(x,u_k) - g(x,\bar{u})\big| dx.
\end{split}
\end{equation}
We note that the boundedness of $\Omega_R$ together with $(G_0)$--$(G_1)$
and \eqref{bddX}, \eqref{limqo2}, Lemma \ref{bddseq1}, allows us to apply
Lebesgue dominated convergence theorem so 
\[
\int_{\Omega_R}\big|g(x,u_k) - g(x,\bar{u})\big| dx \to 0.
\]
On the other hand, again for the boundedness of $\Omega_R$,
by using the arguments in the proof of \cite[Proposition 3.1]{CP2}
from \eqref{bddX}, \eqref{limqo2}, \eqref{strongwpr} and Lemma \ref{bddseq1}
it follows that
\[ 
 \int_{\Omega_R}\big|a(x,u_k,\nabla u_k) - a(x,\bar{u},\nabla \bar{u})\big| |\nabla\varphi| dx\ \to\ 0,
\quad \int_{\Omega_R}\big|A_t(x,u_k,\nabla u_k) - A_t(x,\bar{u},\nabla \bar{u})\big| dx\ \to\ 0.
\]
Thus, summing up, \eqref{crit1} implies \eqref{crit0}.
\end{proof}

\begin{remark}\label{nego}
From Corollary \ref{coroR} we can choose a sequence $(u_k^-)_k$ 
which satisfies all the previous conditions up to \eqref{piu} 
replaced by 
\[
u^-_k(x) \le 0\quad \hbox{for a.e. $x\in \Omega$, for all $k \in \N$.}
\]
Then, arguments, which are similar to those ones developed 
in this section for the proof of Theorem \ref{mainthm},
apply and a critical point $\bar{u}^-$ of $\J$ in $X$ can be found 
so that 
$\bar{u}^-(x) \le 0$ for a.e. $x\in \Omega$.
\end{remark}


\section{Proof of Theorem \ref{mainthm2}}
\label{mainproof2}

Throughout all this section,
we assume that the hypotheses of
Theorem \ref{mainthm2} hold. 
Since this means that stronger assumptions are required 
for $A(x,t,\xi)$ and $g(x,t)$ with respect to those ones in Theorem \ref{mainthm},
all the arguments in Sections \ref{mainbounded} and \ref{mainproof} apply,
and we can consider the sequence $(u_k)_k$, introduced in Section \ref{mainproof},
which satisfies properties $(P_1)$--$(P_3)$, \eqref{cresce}--\eqref{p13} and all the related 
lemmas and propositions.  
Furthermore, Lemma \ref{liminf} can be strengthened as follows.

\begin{lemma}\label{liminfq}
We have that
\[
\liminf_{k\to +\infty} \int_{\Omega} |u_k|^q dx \ >\ 0.
\]
\end{lemma}

\begin{proof}
Arguing by contradiction, assume that 
\[
\liminf_{k\to +\infty} \int_{\Omega} |u_k|^q dx \ =\ 0.
\]
Then, up to subsequences, we have that 
\begin{equation}\label{liminf1q}
\int_{\Omega} |u_k|^q dx \ \to \ 0
\end{equation}
and also 
\begin{equation}\label{liminf2q}
u_k \ \to \ 0 \qquad \hbox{a.e. in $\Omega$.}
\end{equation}
So, taking any $\eps >0$, from $(G_3)$, $(G'_2)$ 
which implies \eqref{altogeps} with $\lambda = 0$,
and \eqref{bddX} we have that
\[
0\ \le\  \int_{\Omega} g(x,u_k) u_k\ dx\ \le\ \eps M_0^p + a_\eps \int_{\Omega} |u_k|^q dx
\]
which, together with \eqref{liminf1q}, implies that
\begin{equation}\label{liminf4q}
\int_{\Omega} g(x,u_k) u_k\ dx\ \to\ 0.
\end{equation}
On the other hand, from \eqref{diff}, \eqref{p13}
and \eqref{liminf4q} it results
\[
\int_{\Omega} \Bigl(a(x,u_k,\nabla u_k)\cdot \nabla u_k 
+ A_t(x,u_k,\nabla u_k) u_k \Bigr)dx\
=\ \int_{\Omega} g(x,u_k)u_k dx\ \to\ 0,
\]
while, from $(H'_5)$ and $(H_3)$, we have that
\[
\int_{\Omega} \Bigl(a(x,u_k,\nabla u_k)\cdot \nabla u_k 
+ A_t(x,u_k,\nabla u_k) u_k \Bigr)dx\
\ge \  \alpha_2\ \left(\alpha_1 \int_{\Omega} |\nabla u_k|^p dx 
+ \int_{\Omega} |u_k|^p dx\right).
\]
Thus, 
\begin{equation}\label{liminf6q}
u_k \ \to \ 0 \quad \hbox{strongly in $W^{1,p}(\Omega)$.}
\end{equation}
At last, from \eqref{bddX}, \eqref{liminf2q} and \eqref{liminf6q},
Proposition \ref{smooth1} applies and 
\[
\J(u_k) \ \to\ \J(0) = 0
\]
in contradiction with property $(P_2)$.
\end{proof}

\begin{remark}\label{estendo}
From $(P_1)$, for each $k \ge 1$ the corresponding function $u_k \in X$
can be trivially extended to the whole Euclidean space $\R^N$
just assuming $u_k(x) = 0$ for all $x \in \R^N \setminus \Omega$.
\end{remark}

In order to prove Theorem \ref{mainthm2}, we need 
the following version of the concentration--compactness lemma (cf. \cite{Li})
which has been stated in \cite[Lemma 1.17]{SZ} if $p=2$.

\begin{lemma}\label{szou}
If $R > 0$ exists such that
\[
\sup_{y \in \R^N}\int_{B_R(y)} |u_k|^p dx\ \to\ 0,
\]
then
\[
u_k \to 0\; \hbox{in $L^q(\R^N)$.}
\]
\end{lemma}

\begin{proof}
From \eqref{bddX}, it is enough reasoning as in the proof of \cite[Lemma 1.17]{SZ}
but replacing 2 with $p > 1$ as $p < q < p^*$.
\end{proof}

Now, we are ready to complete the proof of our dichotomy result.

\begin{proof}[Proof of Theorem \ref{mainthm2}]
From the proof of Theorem \ref{mainthm}
we know that $\bar{u} \in X$, 
weak limit of the sequence $(u_k)_k$ in $W^{1,p}(\Omega)$,
is a critical point of $\J(u)$ in $X$ so that
\eqref{strong1}, \eqref{limqo2} and Proposition \ref{limwpr} hold.\\ 
Firstly, we claim that a constant $\bar{\lambda} > 0$ exists
such that 
\begin{equation}\label{sup1}
\sup_{y \in \R^N}\int_{B_1(y)} |u_k|^p dx\ >\ \bar{\lambda}
\quad \hbox{for all $k \ge 1$.}
\end{equation}
Indeed, otherwise, arguing by contradiction,  
\[
\sup_{y \in \R^N}\int_{B_1(y)} |u_k|^p dx\ \to\ 0
\]
up to subsequences, and so, from Lemma \ref{szou},
it has to be $u_k \to 0$ in $L^q(\Omega)$ in contradiction with 
Lemma \ref{liminfq}.\\
Now, from \eqref{sup1} for any $k \ge 1$ a point $y_k \in \R^N$ exists
such that
\begin{equation}\label{sup2}
\int_{B_1(y_k)} |u_k|^p dx\ >\ \bar{\lambda}.
\end{equation}
Then, two cases may occur: either $(y_k)_k$ has a converging subsequence
or 
\begin{equation}\label{inftyk}
|y_k| \to +\infty.
\end{equation} 
If a subsequence $(y_{k_n})_n$ converges, then a radius $R \ge 1$
exists such that $B_1(y_{k_n}) \subset B_R$ for all $n \ge 1$,
then from \eqref{sup2} and Remark \ref{estendo} it follows that
\begin{equation}\label{sup3}
\int_{\Omega_R} |u_{k_n}|^p dx\ >\ \bar{\lambda} \quad \hbox{for all $n \ge 1$.}
\end{equation}
Thus, \eqref{strong} and \eqref{sup3} imply that  
\[
\int_{\Omega_R} |\bar{u}|^p dx\ \ge\ \bar{\lambda}; 
\]
hence, it has to be $\bar{u} \ne 0$ and $\bar{u}$ is a nontrivial positive solution of
problem \eqref{euler}.\\
On the other hand, if \eqref{inftyk} holds, the statement \eqref{dichotomy}
of this dichotomy result is true.
\end{proof}

\begin{remark}
From Corollary \ref{nego}, the same dichotomy result
can be stated if we consider the sequence of 
nontrivial negative functions $(u_k^-)_k$, 
so either problem \eqref{euler} has a nontrivial negative solution 
or such a sequence satisfies \eqref{dichotomy}.
\end{remark}


\end{document}